\documentclass[11pt]{article}

\usepackage{mathptmx}  % Tang added

\usepackage{amsmath}
\usepackage{bm}
\usepackage{graphicx}
\usepackage{multicol}
\usepackage{wrapfig}
\usepackage{hyperref}
\usepackage{breakurl}
\usepackage{tabularx}
\usepackage{setspace}
\usepackage{comment}
\usepackage{color}
\usepackage{subfig}
\usepackage{cite}

\usepackage{xfrac}
\usepackage{soul}
\usepackage{scrextend}
\usepackage[english]{babel}
\usepackage[autostyle]{csquotes}

\usepackage{sectsty}
\sectionfont{\fontsize{12}{15}\selectfont}

\usepackage{caption}
\usepackage{afterpage}
\usepackage{placeins}
\usepackage{float}

\usepackage{enumitem}
\usepackage{blindtext}

\usepackage[compact]{titlesec} 
\titlespacing{\section}{0pt}{10pt}{5pt}
\titlespacing{\subsection}{0pt}{5pt}{3pt}
\titlespacing{\subsubsection}{0pt}{4pt}{2pt}

\usepackage[margin=1in]{geometry}

\expandafter\def\expandafter\normalsize\expandafter{%
    \normalsize
    \setlength\abovedisplayskip{5pt}
    \setlength\belowdisplayskip{5pt}
    \setlength\abovedisplayshortskip{5pt}
    \setlength\belowdisplayshortskip{5pt}
}

\begin{document}

\begin{center}
{\large\bf Convergence Analysis of Waveform Relaxation Method to 

Compute Coupled Advection-Diffusion-Reaction Equations} 

\vspace*{0.3cm}

W. B. Dong, H. S. Tang\footnote{To whom for correspondence: Hasnong Tang, htang@ccny.cuny.edu}

Department of Civil Engineering, City College

The City University of New York, New York 10031, USA

\end{center}

\vspace*{0.0cm}

\begin{center}
\parbox{0.9\hsize}
{\small

{\sf

{\bf ABSTRACT} \  We study the computation of coupled advection-diffusion-reaction equations by the Schwarz waveform relaxation method. The study starts with linear equations, and it analyzes the convergence of the computation with a Dirichlet condition, a Robin condition, and a combination of them as the transmission conditions. Then, an optimized algorithm for the Dirichlet condition is presented to accelerate the convergence, and numerical examples show a substantial speedup in the convergence. Furthermore, the optimized algorithm is extended to the computation of nonlinear equations, including the viscous Burgers equation, and numerical experiments indicate the algorithm may largely remain effective in the speedup of convergence. 

\vspace{0.1cm}

}
}
\end{center}

\vspace{0.3cm} 

\section {\bf Introduction} 

Advection-diffusion-reaction (ADR) equations and their coupling describe advection, diffusion, reaction, and the interaction among them, which underlie these various interesting and complex phenomena. ADR equations have been solved numerically in studying these phenomena, such as a diffusion-reaction equation depicting migration of pollutants in a porous medium, an advection equation for evolution of ocean surface in the ocean, and an ADR equation involved in an electroosmotic process through a micro-reactor \cite {Atis2012, Qu2019, Esmail2020}. Additionally, ADR equations are also encountered in real-world problems, such as biological process in bone regeneration, heat transfer in phase change, and mechanobiological processes in biomass growth \cite {Prokharau2013, Oliveira2016, Sacco2017}. Moreover, ADR equations and their coupling are used as models to understand complicated partial differential equations (PDEs) and their coupling, respectively \cite {Sakamoto2015, Main2018, Manzanero2020}, such as the Navier Stokes equations, the geophysical fluid dynamics equations, and their coupling \cite {Tang2014,Blayo2016}. Therefore, a study on the computation of ADR equations and particularly their coupling is beneficial to simulating of physical problems and designing numerical PDE methods for them. 
 
ADR problems occur when advection, diffusion, and reaction in a subdomain are coupled with those in another subdomain. Computation of these problems falls into domain decomposition methods, and it traces back over 30 years ago, e.g., \cite {Gastaldi1989}. The methods to compute the problems can be divided into two categories. In the first category, Schwarz iteration is adopted between two adjacent time steps, and hereafter it is referred to as the conventional method. The conventional method is widely applied in the computation of practical problems \cite {Tang2014,Blayo2016}. Particularly, when marching from a time step to the next one, Schwarz iteration is made between the solutions in the two subdomains, exchanging the solutions at subdomain interfaces. When the iteration converges, the solutions are achieved at the next step. The computation repeats such processes to find solutions at all time steps. In the computation, the iteration can be handled by methods for elliptic equations (e.g., \cite {Meurant1991}). Studies on the conventional method include preconditioners \cite {Meurant1991}, algorithms \cite {Canuto2019}, error estimate \cite {Boulaaras2018}, and convergence analysis \cite {Dong2022a}. An interesting result is that an optimal transmission algorithm leads to `perfect convergence', referring to convergence within two times of iteration, and convergence speedup largely remains true in nonlinear situations \cite {Dong2022a}. 

The second category is the waveform relaxation method. In this category, Schwarz iteration is made between solutions in the two subdomains simultaneously at all time steps, rather than separately between two adjacent time steps as in the above conventional method. The waveform relaxation method traces back to over two decades ago, e.g., in attempts to solve systems of algebraic equations and differential equations \cite {Lelarasmee1982,Crow1994, Gander:2002}, and it has evolved into a popular method \cite {Gander2008}. Study on this method includes transmission conditions, convergence analysis, optimized transmission algorithms, and computation of equations with discontinuous coefficients \cite {Martin2004, Gander2007a, Gander2007b, Califano2018}. A substantial understanding has been achieved from the past investigations. For instance, it is indicated that transmission conditions have a strong influence on convergence speed \cite {Gander2007a}. The computation with Dirichlet conditions at interfaces is well-posed, and the convergence rate increases with the size of the overlapping region but decreases with diffusion coefficients \cite {Gander2007b}. An obstacle is that the convergence rate gets slower as the grid spacing gets finer \cite {Gastaldi1989, Gander1998}. Computation by a Schwarz waveform relaxation method is studied for systems of ordinary differential equations resulting from {RC} type circuits \cite {Gander2004, AlKhaleel2014}, which seem different but are related to the discretization of coupled advection-diffusion reactions equations as shown in this paper. 

We have analyzed the computation of coupled ADR equations by the conventional method \cite {Dong2022a}, and this paper makes a further study on the computation of the equations by the waveform relaxation method, aiming to further understand the roles of transmission conditions as well as their algorithms. Most earlier investigations deal with the same linear equations in subdomains, and they are primarily conducted at the continuous level. e.g., \cite {Gander2007b,Califano2018}. Different from the existing investigations, this work presents a study on a coupling of ADR equations, which may be the same or different in subdomains. Additionally, it deals with semi-discretization systems of the equations. Dirichlet condition, Robin condition, and a  combination of them are adopted as the transmission conditions. In reference to the techniques on RC-type circuits in \cite  {Gander2004, AlKhaleel2014}, an optimized transmission algorithm is analyzed for the transmission conditions. Moreover, a discussion is made on computation of nonlinear equations, and exploration is made on validity of the optimized algorithm in computation of these equations. 
%It is anticipate the effort of this paper promotes the  development of computational methods for advection-diffusion-reaction equations.

\section {\bf The problem of computation and transmission condition}

\ul {Equation, discretization, and solution} Consider a linear advection-diffusion-reaction (ADR) equation as follows:
\begin{equation}
\label{equation}
u_t + \mu u_x = \theta u_{xx}-\gamma u
\end{equation}
\noindent in which $\mu  $, $\theta $, $\gamma =consts$, $\mu \ne 0$, $\theta $, $\gamma >0$. By central difference in space on a grid $(\cdots, i-1, i, i+1, \cdots )$, the equation can be transferred into a following semi-discretization system:   
\begin{equation}
\label{discret1}
\begin{array}{l}
u_t + \mu  \dfrac {u_{i+1}-u_{i-1}}{2\Delta x} = \theta \dfrac {u_{i+1}-2u_i+u_{i-1}}{\Delta x^2}-\gamma u_i,  \ \ 
\end{array}
\end{equation}

\noindent where $\Delta x$ is grid spacing. After the Laplace transform (for simplicity, consider $u_i=0$, at $t=0$), Eq. (\ref {discret1}) can be expressed as the following linear system:
\begin{equation}
\label{Laplace}
a\hat {u}_{i-1}+(b-s)\hat {u}_i+c\hat {u}_{i+1}=0
\end{equation}
\noindent in which the hat  `$\string^$' stands for the Laplace transform, and $s=\sigma +i\omega $, being a complex number, with $\sigma \ge 0$ is considered. Moreover, 
\begin{equation}
\label{abc}
a=\frac {\mu  }{2\Delta x}+\frac {\theta }{{\Delta x}^2}, \ \
b=-\frac {2\theta }{{\Delta x}^2}-\gamma , \ \ 
c=-\frac {\mu  }{2\Delta x}+\frac {\theta }{{\Delta x}^2}
\end{equation}
\noindent In view of the values for $\mu  $, $\theta $, and $\gamma $, it is seen that $b<0$. In this study, we consider $c\ne 0$ only. 

The linear system (\ref {Laplace}) has a solution in the following form: 
\begin{equation}
\label{gen_solu}
\hat {u}=A{r_-}^i+B{r_+}^i 
\end{equation}

\noindent as long as 
\begin{equation}
\label {r_-+}
|r_-|< 1, \qquad |r_+|> 1 
\end{equation}

\noindent in which $A$, $B$ $=const$, and 
\begin{equation}
\label{roots1}
r_\pm =\frac {s-b\pm \sqrt {(s-b)^2-4ac}}{2c}
\end{equation}

\noindent which are the roots of the characteristic equation 
\begin{equation*}
\label{characteristic_eq}
cr^2+(b-s)r+a=0. 
\end{equation*}

\vspace {0.3 cm}

{\textsc{proposition 1.1}} The roots in (\ref{roots1}) are analytical. Furthermore, $|r_+|>1$, and, under condition
\begin{equation}
\Big |i\omega -b+ \sqrt {(i\omega -b)^2-4ac}\Big |>2a
\label {cond_r_-}
\end{equation} 

\noindent $|r_-|< 1$ holds. 

\vspace {0.3 cm}

The proof of the proposition is given in the appendix, and situations in which condition (\ref {cond_r_-}) holds are described there. This study considers the situations under condition (\ref {r_-+}).   

\begin{figure}[h]
\vspace {-0. cm}
\center
\includegraphics[width=0.4\linewidth]{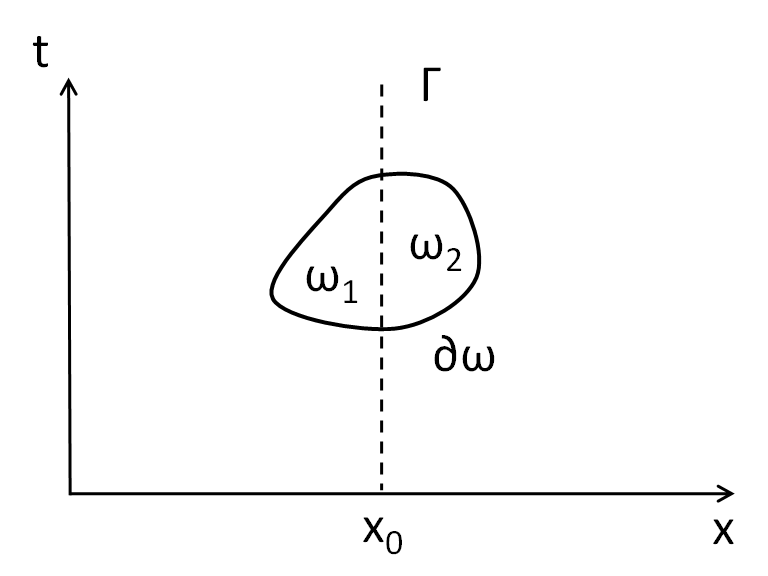}
\vspace {-0. cm}
\caption{\small Domain division.}
\label{fig1}
\vspace {-0. cm}
\end{figure}

\ul {Problem of computation and transmission condition} Let the plane $x-t$ be divided by an interface $\Gamma $, which is located at $x=x_0$, into two parts, as shown in Fig. \ref {fig1}. Then, we consider a problem consisting of two coupled initial value problems of Eq. (\ref {equation}) as follows:
\begin{equation}
\label{ivp0}
\left\{\begin{array}{ll}
v_t + \mu _1v_x = \theta _1v_{xx}-\gamma _1v,  \ & t>0 \\
v=f(x),  \ & t=0 \\
v=p_1(w), \ & x=x_0 
\end{array} \right.  \qquad 
\left\{\begin{array}{ll}
w_t + \mu _2w_x = \theta _2w_{xx}-\gamma_2w, \ & t>0  \\
w=f(x), \ & t=0 \\ 
w=p_2(v), \ & x=x_0 
\end{array} \right. 
\end{equation}

\noindent where $\mu _l$, $\theta _l$, $\gamma _l = consts$, $\mu \ne 0$, $\theta _l$, $\gamma _l$ $> 0$ ($l=1,2$). 

$p_l$ are operators for the transmission conditions. 
%\begin{wrapfigure}{r}{0.43\textwidth} % l, r, c; this letters stand for left, right, centre,
%\end{wrapfigure}

An transmission condition is the continuity of the solutions of the two equations at the interface, that is
\begin{eqnarray}
v=w, \ \ x=x_0
\label{intf_cond1}
\end{eqnarray}

\noindent Such condition is commonly used in computation of practical problems. Let us consider another type of transmission condition. Suppose $\omega $ is an any region across the interface, and $\omega =\omega _1\cup \omega _2 $, with $\omega _1$ and $\omega _2 $ falling in left and right side of $x=x_0 $, respectively (Fig. \ref {fig1}). Let $u=v\cup w$, $\mu  =\mu  _1\cup \mu  _2$, and $\theta =\theta _1\cup \theta _2$, we consider the weak form of the equations in (\ref {ivp0}) as follows, 
\begin{eqnarray*}
\int _\omega u\phi _t+(\mu  u-\theta u_x)\phi _x+\gamma u\phi )dxdt
%+\int _{-\infty }^{\infty} (u\phi )\vert_{t=0}dx
=0  
\label{weak_solu1}
\end{eqnarray*}

\noindent in which $\forall \phi \in C_0^\infty  $. By integral by part, one has 
\begin{eqnarray*}
\begin{array}{ll}
&\int _{\omega _1}(-v_t-\mu  _1v_x+\theta v_{xx}-\gamma _1v)\phi dxdt +\int _{\omega _2}(-w_t-\mu  _2 w_x+\theta _2w_{xx}-\gamma _2w)\phi dxdt \\
&+\int _\Gamma (\mu  _1v-\theta v_x)\vert _{\Gamma _-}\phi dt -\int _\Gamma (\mu  _2w-\theta _2 w_x)\vert _{\Gamma _+}\phi dt =0
\end{array}
\label{weak_solu2}
\end{eqnarray*} 

\noindent The parentheses in the first two terms become zero because of the equations in (\ref {ivp0}). Since $\omega $ is arbitrarily selected, above leads the following Robin condition at the interface: 
\begin{eqnarray}
g(v)|_{\Gamma _-}=h(w)|_{\Gamma _+}  
\label{intf_cond2}
\end{eqnarray}

\noindent where $g(v)=\mu  _1 v-\theta _1 v_x$, $h(w)=\mu  _2 w-\theta _2 w_x$. $g$ and $h$ are the fluxes of the equations in (\ref {ivp0}), and thus (\ref {intf_cond2}) requires that the fluxes are continuous across the interface. Note that condition (\ref {intf_cond2}) does not guarantee condition \ref {intf_cond1}, or, it permits a discontinuous solution across the interface. 

We consider two subdomains, one on the left: $x<x_2$, one on the right: $x>x_1$, with $x_1<x_2$. The two subdomains overlap within $x_2<x<x_1$, and their interface locations are at $x=x_1$, $x_2$. In this work, transmission conditions based on above two transmission conditions and their combination will be considered. The first set of transmission conditions is 
\begin{equation}
\label{D_cond}
v=w, \  x=x_2; \qquad w=v, \ x=x_1 
\end{equation}

\noindent which imposes continuity of solutions at the two interfaces via Dirichlet conditions. Hereafter, this set of conditions will be referred to as a Dirichlet condition. The second set of transmission conditions is  
\begin{equation}
\label{Robin_cond}
\begin{array}{ll}
g(v)=h(w),   \  x=x_2; \qquad  h(w)=g(v), \ x=x_1
\end{array}
\end{equation}

\noindent which imposes continuity of fluxes at the interfaces. Such set of transmission conditions will be referred as to a Robin condition. The third set of transmission conditions is
\begin{equation}
\label{comb_cond}
v=w,   \  x=x_2; \qquad   h(w)=g(v), \ x=x_1
\end{equation}

\noindent which is a combination of the Dirichlet and Robin conditions, or, a combined Dirichlet and Robin condition, and it will be referred to as the combined condition. 

\section {Computation of coupled equations}

\ul {Computation and Schwarz iteration} Now, let us consider computation of coupled problems (\ref {ivp0}) via Schwarz iteration as follows:   
\begin{equation}
\label{ivp1}
\left\{\begin{array}{ll}
v_t^{k+1} + \mu _1v_x^{k+1} = \theta _1v_{xx}^{k+1}-\gamma _1v^{k+1},  \ & t \in (0,T] \\
v^{k+1}=f(x),  \ & t=0 \\
v^{k+1}=p_1(w^k), \ & x=x_2 
\end{array} \right.  \quad  
\left\{\begin{array}{ll}
w_t^{k+1} + \mu _2w_x^{k+1} = \theta _2w_{xx}^{k+1}-\gamma_2w^{k+1}, \ & t \in (0,T]  \\
w^{k+1}=g(x), \ & t=0 \\ 
w^{k+1}=p_2(v^k), \ & x=x_1 
\end{array} \right. 
\end{equation}
\noindent The computation is made on two grids, one on the left indexed with ($i= ... -2, -1, 0$), and one on the right indexed with ($i=0, 1, 2, ... $), see Fig. \ref {fig2}. Utilizing discretization in (\ref {discret1}), the computation is implemented via the following Schwarz iteration ($k\ge 0$, $t \in [0,T]$): 
\begin{figure}[h] 
  \vspace{-0.cm}
  \centering
  {\includegraphics[width=0.48\textwidth]{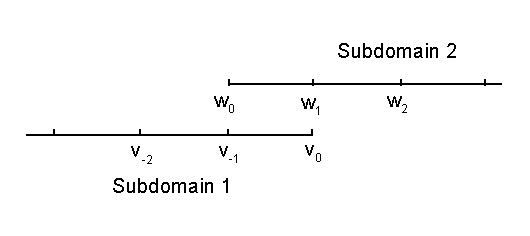}} \ \ \ 
  \vspace{-8pt}
  \caption{\small Two grids overlap at two nodes.}
  \label{fig2}
  \vspace{-0pt}
\end{figure}

\begin{equation}
\label{dis_ivp1}
\left\{\begin{array}{ll}
{v_t}_i^{k+1}=a_1v_{i-1}^{k+1}+b_1v_i^{k+1}+c_1v_{i+1}^{k+1},  \ \ i\le -1 \\
v_0^{k+1}=\tilde {p}_1(w_0^k,w_1^k, ...)  \\
\end{array}\right . \ \ 
\left\{\begin{array}{ll}
{w_t}_i^{k+1}=a_2w_{i-1}^{k+1}+b_2w_i^{k+1}+c_2w_{i+1}^{k+1},  \ \ i\ge 1 \\
w_0^{k+1}=\tilde {p}_2(v_0^k,v_{-1}^k, ...)\\ 
\end{array}\right . 
\end{equation}

\noindent where  $ a_l$, $b_l$, and $c_l$ ($l=1,2$) are defined as in (\ref {abc}), and $\tilde {p}_l$ ($l=1,2$) are discretization of the operators for transmission conditions. Eq. (\ref {dis_ivp1}) can be expressed as two linear systems: 

\begin{equation*}
\label{int_Sch}
\begin{array}{ll}
\begin{bmatrix}
\vdots                   \\
\vdots                   \\
{v_t}_{-2}^{k+1}         \\
{v_t}_{-1}^{k+1}  
\end{bmatrix}
=
\left [
\begin{array}{rrrrrrrrrrrrrrrrrrr}
           &         &\ddots   & \vdots    & \vdots    & \vdots  \\ 
           & \cdots  &   a_1   &   b_1     & c_1       & 0          \\
           & \cdots  &   0     &   a_1     & b_1       & c_1         \\
           &  \cdots  &   0     &   0       & a_1       & b_1 
\end{array} 
\right ]
\begin{bmatrix}
\vdots   \\
\vdots   \\
v_{-2}^{k+1}         \\
v_{-1}^{k+1}  
\end{bmatrix}
+
\begin{bmatrix}
   \vdots   \\
   \vdots   \\
   0        \\    
 c_1v_0^{k+1}
\end{bmatrix}  \\[40pt]

\begin{bmatrix}
{w_t}_{-1}^{k+1}     \\
{w_t}_{-2}^{k+1}     \\
\vdots   \\ 
\vdots \\
\end{bmatrix} 
=
\left [
\begin{array}{rrrrrrrrrrrrrrrrrrr}
 &   b_2       &   c_2     & 0         & 0          & \cdots   \\
 &   a_2       &   b_2     & c_2       & 0          & \cdots   \\
 &   0         &   a_2     & b_2       & c_2        & \cdots   \\
 &   \vdots    & \vdots    & \vdots    & \ddots     &          \\ 

\end{array} 
\right ]
\begin{bmatrix}
w_1^{k+1}     \\
w_2^{k+1}     \\
\vdots   \\
\vdots \\
\end{bmatrix}
+
\begin{bmatrix}
 a_2w_0^{k+1}  \\  
   0      \\ 
   \vdots   \\
   \vdots \\
\end{bmatrix}
\end{array}
\end{equation*}

\noindent Making a Laplace transform on above iteration, one has

\begin{equation}
\label{trans_int_Sch}
\begin{array}{ll}
\left [
\begin{array}{rrrrrrrrrrrrrrrrrrr}
          & \ddots  &\vdots   & \vdots    & \vdots    & \vdots   \\ 
          & \cdots  &   a_1   &   b_1-s   & c_1       & 0        \\
          & \cdots  &   0     &   a_1     & b_1-s     & c_1      \\
          & \cdots  &   0     &   0       & a_1       & b_1-s 
\end{array} 
\right ]
\begin{bmatrix}
\vdots  \\
\hat {v}_{-3}^{k+1}  \\
\hat {v}_{-2}^{k+1}  \\
\hat {v}_{-1}^{k+1}  
\end{bmatrix}
=
\begin{bmatrix}
   \vdots   \\
   \vdots \\
   0      \\
-c_1\hat {v}^{k+1}_0
\end{bmatrix}  \\[30pt]
\left [
\begin{array}{rrrrrrrrrrrrrrrrrrr}
 b_2-s     &   c_2     & 0         & 0          & \cdots    &    \\
 a_2       &   b_2-s   & c_2       & 0          & \cdots    &    \\
 0         &   a_2     & b_2-s     & c_2        & \cdots    &    \\
 \vdots    & \vdots    & \vdots    & \vdots     & \ddots    &    \\ 

\end{array} 
\right ]
\begin{bmatrix}
\hat {w}_1^{k+1}     \\
\hat {w}_2^{k+1}     \\
\hat {w}_3^{k+1}     \\
\vdots    \\
\end{bmatrix}
=
\begin{bmatrix}
-a_2\hat {w}^{k+1}_0  \\  
   0      \\ 
   \vdots \\
   \vdots \\
\end{bmatrix}
\end{array}
\end{equation}

\noindent The solution to the first system is in form of (\ref {gen_solu}), and, in view that $|r_-|<1$, the coefficient $A$ becomes zero. Similarly, the solution to the second system is also in form of (\ref {gen_solu}), with $B$ to be zero. As a result, these yield to 
\begin{equation}
\label{int_Sch1_solu}
\begin{array}{l}
\hat {v}_i^{k+1}=B^{k+1}{r_+}^i,  \ \  i=0, -1, -2, ... \\
\hat {w}_i^{k+1}=A^{k+1}{r_-}^i,  \ \  i=0,  1,  2, ...
\end{array}
\end{equation}

\noindent Here and hereafter $r_+$ and $r_-$ are evaluated by $(a_1,b_1,c_1)$ and $(a_2,b_2,c_2)$, respectively. Based on the above discussion, the convergence in computation of (\ref {dis_ivp1}) in association with different transmission conditions are discussed in the following content. 

\ul {Dirichlet condition} Applying of Dirichlet condition (\ref {D_cond}) at the interfaces of the two grids when they overlap at two nodes (Fig. \ref {fig2}), the transmission condition in (\ref {dis_ivp1}) becomes 
\begin{equation}
\label{numerical_Dirichlet_cond}
v_0^{k+1}=w_1^k, \ \ \ w_0^{k+1}=v_{-1}^k
\end{equation}

\noindent which, after Laplace transform, becomes
\begin{equation*}
\label{transformed_numerical_Dirichlet_cond}
\hat {v}_0^{k+1}=\hat {w_1}^k, \ \ \ \hat {w}_0^{k+1}=\hat {v}_{-1}^k
\end{equation*}

\noindent Now, with the aid of above, the last equation in the first system in (\ref {trans_int_Sch}) becomes

\begin{equation*}
a_1\hat {v}_{-2}^{k+1}+(b_1-s)\hat {v}_{-1}^{k+1}+c_1\hat {w}_1^k=0
\end{equation*}

\noindent which, together with (\ref {int_Sch1_solu}), leads to 

\begin{equation*}
B^{k+1}=-\dfrac {c_1r_-}{a_1{r_+}^{-2}+(b_1-s){r_+}^{-1}}A^k 
\end{equation*}

\noindent Similarly, one has

\begin{equation*}
A^{k+1}=-\dfrac {a_2{r_+}^{-1}}{c_2{r_-}^2+(b_2-s)r_-}B^k
\end{equation*}

\noindent From the above equations, it is readily derived that
\begin{equation}
\label{Laplace_iter_conv}
\begin{array}{l}
\hat {v}_i^{k+1}=\rho \hat {v}_i^{k-1}, \ \ \ i=0, -1, -2, ... \\
\hat {w}_i^{k+1}=\rho \hat {w}_i^{k-1}, \ \ \ i=0,  1,  2, ... \\
\end{array}
\end{equation}

\noindent in which
\begin{equation*}
\rho =\dfrac {a_2c_1}{(c_2r_-+b_2-s)(a_1{r_+}^{-1}+b_1-s)}
\end{equation*}

\noindent Here, $\rho $ is the contraction factor of the computation, which reflects the convergence speed. Furthermore, in view of (\ref {roots1}) and (\ref {roots4}), one has 
\begin{equation}
\label {work1}
\begin{array}{ll}
&a{r_+}^{-1}+b-s \\ 
=&c r_-+b-s \\
=&c\dfrac {s-b- \sqrt {(s-b)^2-4ac}}{2c}+b-s \\
=&-c r_+\\
=&-a{r_-}^{-1}
\end{array}
\end{equation}

\noindent the contraction factor becomes
\begin{equation}
\label{Rho_Dirichlet}
\begin{array}{ll}
\rho &=\dfrac {a_2c_1}{(-a_2{r_-}^{-1})(-c_1r_+)}  \\
     &=\dfrac {r_-}{r_+} 
\end{array}
\end{equation}

\noindent In view of (\ref {r_-+}) and all above discussion, we arrive at the following conclusion. 

\vspace {0.3 cm}

\textsc {Theorem 3.1} When (\ref {dis_ivp1}) is computed via Dirichlet transmission condition (\ref {numerical_Dirichlet_cond}), the contraction rate is $\rho =r_-/r_+$. Moreover, $|\rho |<1$, that is, the computation converges.  

\vspace {0.3 cm}

\textsc {Remark 3.1} Above theorem holds for different scenarios, including coupling of two same equations, two equations of a same type but with different coefficients, and two equations of different types (e.g., an advection equation and an advection-diffusion equation). Theorem 3.1 indicates that computations of all of these scenarios in conjunction of transmission condition  (\ref {numerical_Dirichlet_cond}) will converge. 

\vspace {0.3 cm}

\textsc {Remark 3.2} The conclusions in Theorem 3.1 has been reported in a study on the RC-type circuits with $a_1=a_2=c_1=c_2$, $b_1=b_2$ \cite {Gander:2004}, and thus the theorem may be considered as an extension of the previous work.  

\vspace {0.3 cm}

\ul {Robin condition} Now, consider the computation of (\ref {dis_ivp1}) in association of Robin condition (\ref {Robin_cond}). When the two grids overlap at two grid nodes (Fig. \ref {fig2}), the Robin condition becomes 
\begin{equation}
\label{numerical_Robin1_cond}
\begin{array}{l}
\bar {g}_{-1/2}^{k+1}=\bar {h}_{1/2}^k, \ \bar {h}_{1/2}^{k+1}=\bar {g}_{-1/2}^k  
\end{array}
\end{equation}

\noindent where $\bar g$ and $\bar h$ are an approximation of flux $g$ and $h$, respectively, in condition (\ref {Robin_cond}), and they are numerical fluxes:  
\begin{equation*}
\label{numer_flux}
\begin{array}{l}
\bar {g}_{-1/2}=\dfrac {\mu  _1(v_{-1}+v_0)}{2}-\dfrac {\theta _1(v_0-v_{-1})}{2} \\
\bar {h}_{1/2}=\dfrac {\mu  _2(w_0+w_1)}{2}-\dfrac {\theta _2(w_1-w_0)}{2}
\end{array}
\end{equation*}

\noindent The transmission algorithm leads to
\begin{equation*}
\begin{array}{l}
v_0^{k+1}=\dfrac {a_1}{c_1}{v_-}_1^{k+1}-\dfrac {a_2}{c_1}w_0^k+\dfrac {c_2}{c_1}w_1^k \\
w_0^{k+1}=\dfrac {c_2}{a_2}w_1^{k+1}+\dfrac {a_1}{a_2}v_{-1}^k-\dfrac {c_1}{a_2}v_0^k 
\end{array}
\end{equation*}

\noindent Then, with the aid of above, the last equation in the first system in (\ref {trans_int_Sch}) becomes
\begin{equation*}
a_1\hat {v}_{-2}^{k+1}+(b_1-s)\hat {v}_{-1}^{k+1}
+c_1\left ( \dfrac {a_1}{c_1} \hat {v}_{-1}^{k+1}-\dfrac {a_2}{c_1}\hat {w}_0^k+\dfrac {c_2}{c_1}\hat {w}_1^k\right )=0
\end{equation*}

\noindent Using solution (\ref {int_Sch1_solu}), the above equation leads to 
\begin{equation*}
B^{k+1}=\dfrac {a_2-c_2r_-}{a_1r_+^{-2}+(b_1-s)r_+^{-1}+a_1r_+^{-1}}A^k,  
\end{equation*}

\noindent Similarly, one has
\begin{equation*}
A^{k+1}=\dfrac {c_1-a_1{r_+}^{-1}}{c_2r_-+(b_2-s)r_-+c_2r_-^2}B^k, 
\end{equation*}

\noindent Again, according to (\ref {int_Sch1_solu}), the last two equations give rise to 
\begin{equation*}
\begin{array}{l}
\rho =\dfrac {a_2-c_2{r_-}}{a_1{r_+}^{-2}+(b_1-s){r_+}^{-1}+a_1{r_+}^{-1}} \cdot 
      \dfrac {c_1-a_1{r_+}^{-1}}{c_2{r_-}+(b_2-s){r_-}+c_2{r_-}^2}
\end{array}
\end{equation*}

\noindent Then, with the aid of (\ref {work1}) and (\ref {roots4}), it is derived that  
\begin{equation}
\label{rho_Robin1}
\begin{array}{ll}
\rho &= \dfrac {a_2-c_2r_-}{-c_1r_+{r_+}^{-1}+a_1{r_+}^{-1}}\cdot
         \dfrac {c_1-a_1{r_+}^{-1}}{c_2r_--a_2{r_-}^{-1}r_-} \\
     &=1
\end{array}
\end{equation}

\noindent Therefore, Schwarz iteration has no convergence when transmission algorithm (\ref {numerical_Robin1_cond}) is adopted. 

Now, consider a slight modification to (\ref {numerical_Robin1_cond}). Instead of two nodes, let the two grids overlap at three nodes (i.e., $v_{-2}$ and $w_0$ are at a same location, and $v_0$ and $w_2$ are at the same location), and apply the Robin condition as follows, 
\begin{equation}
\label{numerical_Robin2_cond}
\begin{array}{l}
\hat {f}_{-1/2}^{k+1}=\hat {g}_{3/2}^k, \ \ \hat {g}_{1/2}^{k+1}=\hat {f}_{-3/2}^k
\end{array}
\end{equation}

\noindent which leads to 
\begin{equation*}
\begin{array}{l}
v_0^{k+1}=\dfrac {a_1}{c_1}v_{-1}^{k+1}-\dfrac {a_2}{c_1}w_1^k+\dfrac {c_2}{c_1}w_2^k \\
w_0^{k+1}=\dfrac {c_2}{a_2}w_{-1}^{k+1}-\dfrac {c_1}{a_2}v_{-1}^k+\dfrac {a_1}{a_2}v_{-2}^k 
\end{array}
\end{equation*}

\noindent By a similar derivation, one has 
\begin{equation}
\label {Rho_Robin2}
\begin{array}{ll}
\rho &=\left (\dfrac {a_2-c_2r_-}{a_1{r_+}^{-2}+(b_1-s){r_+}^{-1}+a_1{r_+}^{-1}} \cdot 
      \dfrac {c_1-a_1{r_+}^{-1}}{c_2{r_-}+(b_2-s)r_-+c_2{r_-}^2}\right ) \cdot
      \dfrac {r_-}{r_+}\\
     &=\dfrac {r_-}{r_+}  
\end{array}
\end{equation} 

\noindent In above, the terms in the parentheses actually comprise the previous contraction factor with the two grid overlap overlap at two nodes, whose value is unit in magnitude. As a result, the contraction factor is the same to that for the Dirichlet condition (\ref {numerical_Dirichlet_cond}) is used, and, because of (\ref {r_-+}), the computation will converge. The above is summarized as follows. 

\vspace {0.3 cm}

{\textsc{Theorem 3.2}} When (\ref {dis_ivp1}) is computed via the Robin condition (\ref {numerical_Robin1_cond}),  $\rho =1$, and thus the computation has no convergence. However, when the Robin condition (\ref {numerical_Robin2_cond}) is adopted,  $\rho =r_-/r_+$, and $|\rho |<1$, that is, the computation converges. 

\vspace {0.3 cm} 

\ul {Combined condition} Now, let us apply the combined condition, (\ref {comb_cond}). Particularly, consider that the two grids overlap at two nodes, see Fig. \ref {fig2}, and the combined condition is implemented as 
\begin{equation}
\label{numerical_comb2_cond}
\begin{array}{l}
v_0^{k+1}=w_1^k, \ \ \ \hat {h}_{1/2}^{k+1}=\hat {g}_{-1/2}^k  
\end{array}
\end{equation}

%\begin{equation}
%\label{numerical_comb2_cond}
%\begin{array}{l}
%v_0^{k+1}=w_1^k, \ \ w_0^{k+1}=\dfrac {c_2}{a_2}w_1^{k+1}-\dfrac {c_1}{a_2}v_0^k+\dfrac {a_1}{a_2}v_{-1}^k 
%\end{array}
%\end{equation}

\noindent By steps similar to those of the above, it is derived that the contraction factor becomes
\begin{equation*}
\begin{array}{l}
\rho =\dfrac {c_1r_-}{a_1{r_+}^{-2}+(b_1-s){r_+}^{-1}} \cdot
      \dfrac {c_1-a_1{r_+}^{-1}}{c_2r_-+(b_2-s)r_-+c_2{r_-}^2}  
\end{array}
\end{equation*}

\noindent With the aid of (\ref {work1}), it can be derived that 
\begin{equation}
\label{rho_comb2}
\begin{array}{l}
\rho =\dfrac {c_1-a_1{r_+}^{-1}}{c_2-a_2{r_-}^{-1}}
\end{array}
\end{equation}

\noindent It is seen that the combination of the two types of transmission algorithms does not necessarily lead to convergence. From the above, a condition for the convergence in this situation is  
\begin{equation}
\label {stability_comb2_cond}
\begin{array}{l}
\left | \dfrac {a_1-c_1r_+}{a_2-c_2r_-}\right |<\left | \dfrac {r_+}{r_-}\right |
\end{array}
\end{equation}

Now we consider overlapping of three grid nodes, and apply the combined  condition  
\begin{equation}
\label{numerical_comb3_cond}
\begin{array}{l}
v_0^{k+1}=w_2^k, \ \ \ \hat {g}_{1/2}^{k+1}=\hat {h}_{-3/2}^k 
\end{array}
\end{equation}

\noindent In this situation, it is derived that 
\begin{equation}
\label{rho_comb3}
\begin{array}{l}
\rho =\dfrac {c_1-a_1{r_+}^{-1}}{c_2-a_2{r_-}^{-1}}\cdot \dfrac {r_-}{r_+}
\end{array}
\end{equation}

\noindent In above, the first fractional term is the contraction factor for the computation with two overlapping nodes, and the second fractional term is less than 1 in magnitude. Therefore, interestingly, this shows that adding an overlapping node leads to a faster convergence if the computation converges. A convergence condition in this situation reads as
\begin{equation}
\label{stability_comb3_cond}
\left | \dfrac {a_1-c_1r_+}{a_2-c_2r_-}\right |<\left | \dfrac {r_+}{r_-}\right |^2
\end{equation}

All above discussions are summarized as follows. 

\vspace {0.3 cm}

\textsc{Theorem 3.3} When (\ref {dis_ivp1}) is computed via the combined  condition (\ref {numerical_comb2_cond}), a convergence condition is (\ref {stability_comb2_cond}). While it is computed via the combined condition (\ref {numerical_comb3_cond}), a convergence condition is (\ref {stability_comb3_cond}), and, the convergence speed is faster than via the former condition if the computation converges.  

\vspace {0.3 cm}

\begin{figure}[!ht] 
  \vspace{-0.2 cm}
  \centering
  \subfloat[][]
  {\label{fig3a}\includegraphics[width=0.32\textwidth]{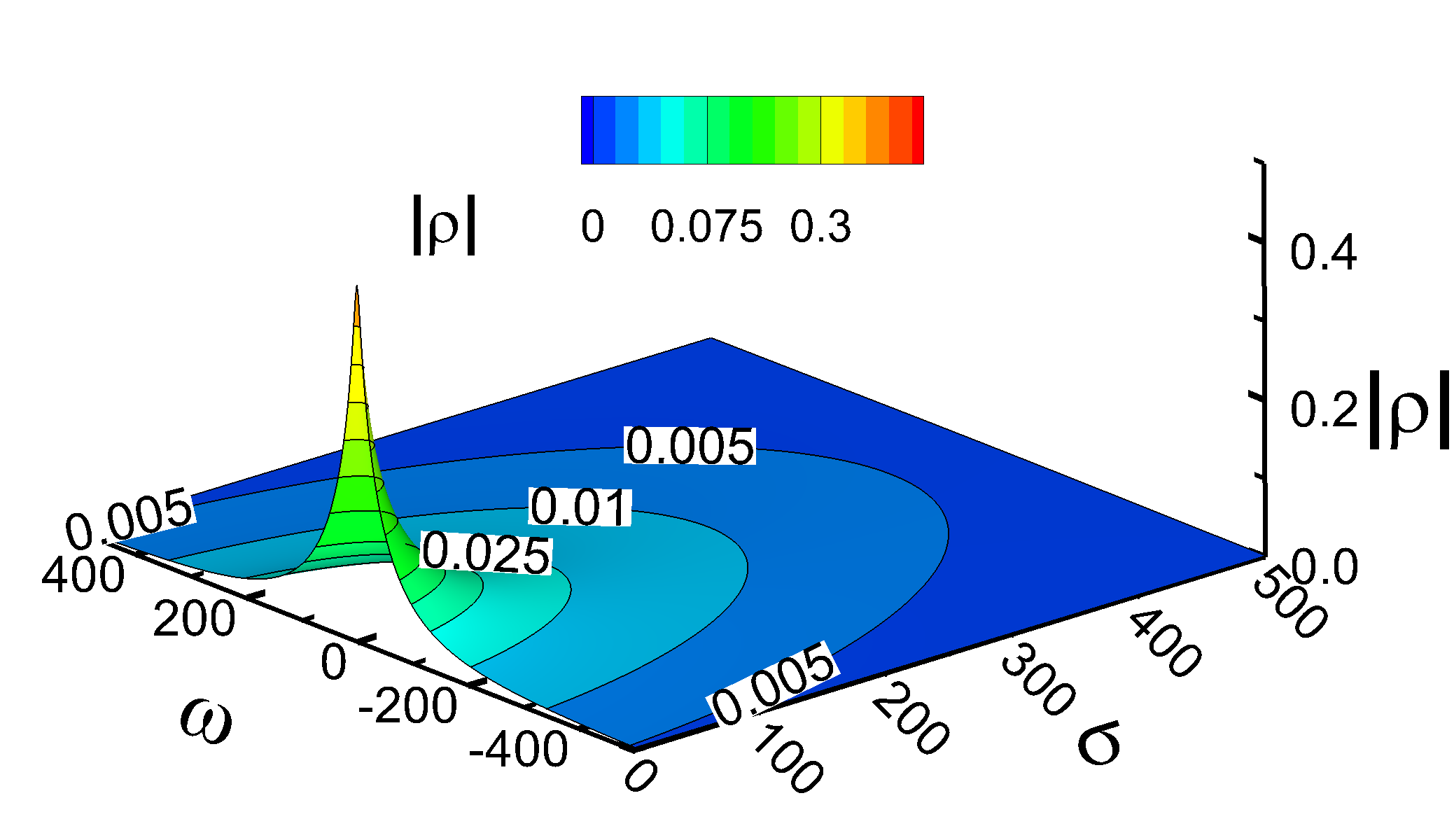}} \ 
  \subfloat[][]
  {\label{fig3b}\includegraphics[width=0.32\textwidth]{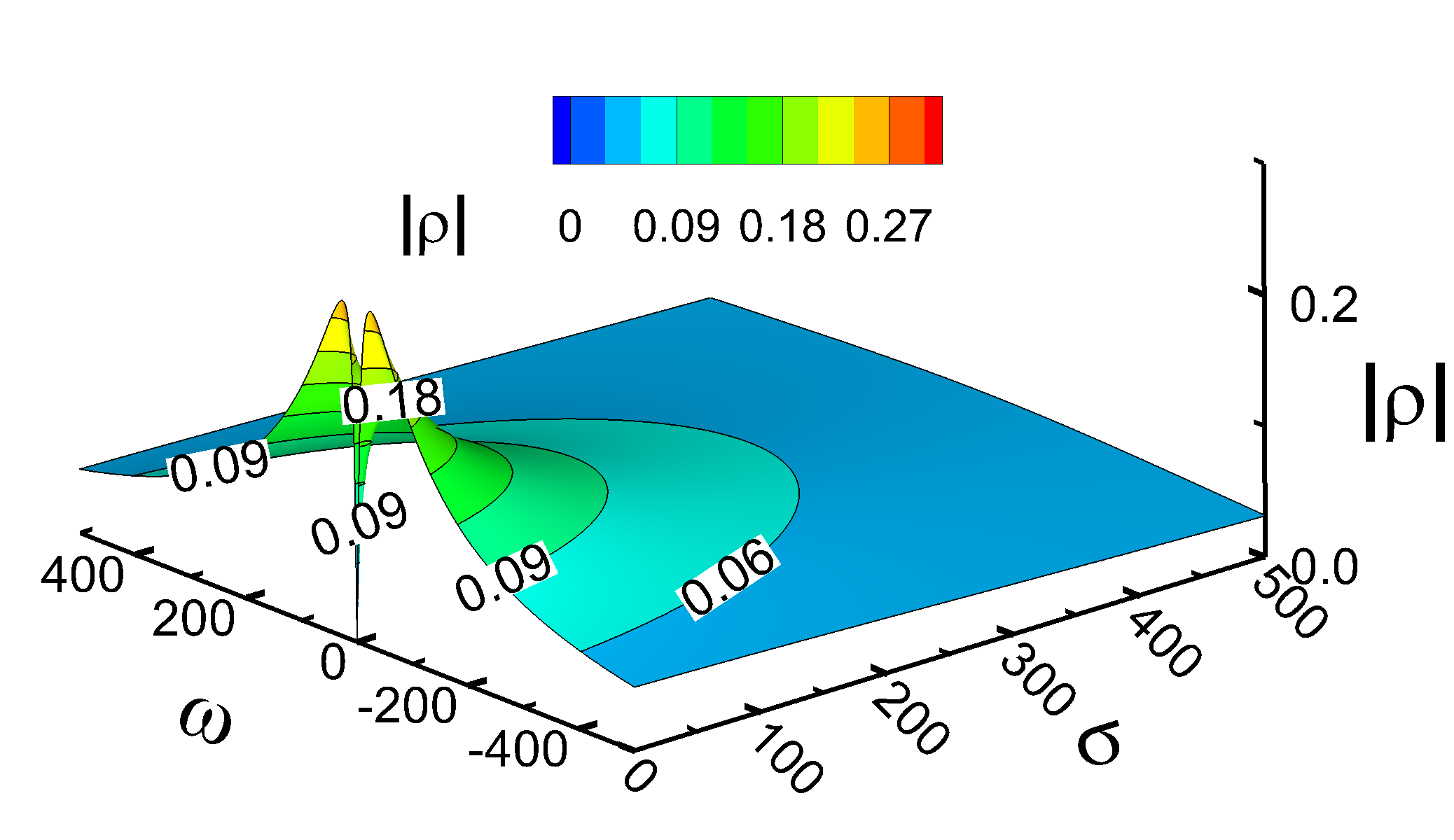}} \ 
  \subfloat[][]
  {\label{fig3c}\includegraphics[width=0.32\textwidth]{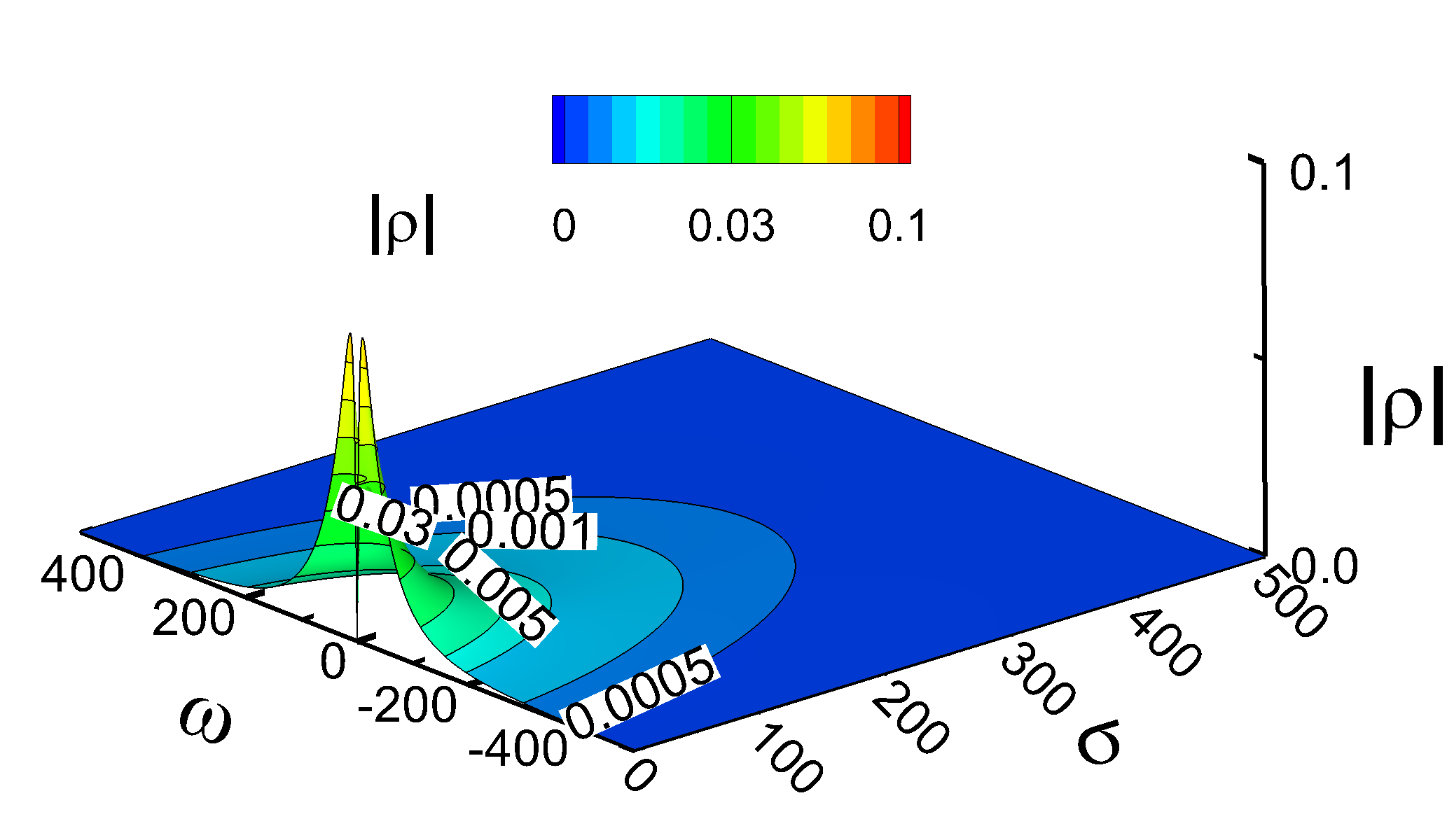}} \\
  \subfloat[][]
  {\label{fig3d}\includegraphics[width=0.32\textwidth]{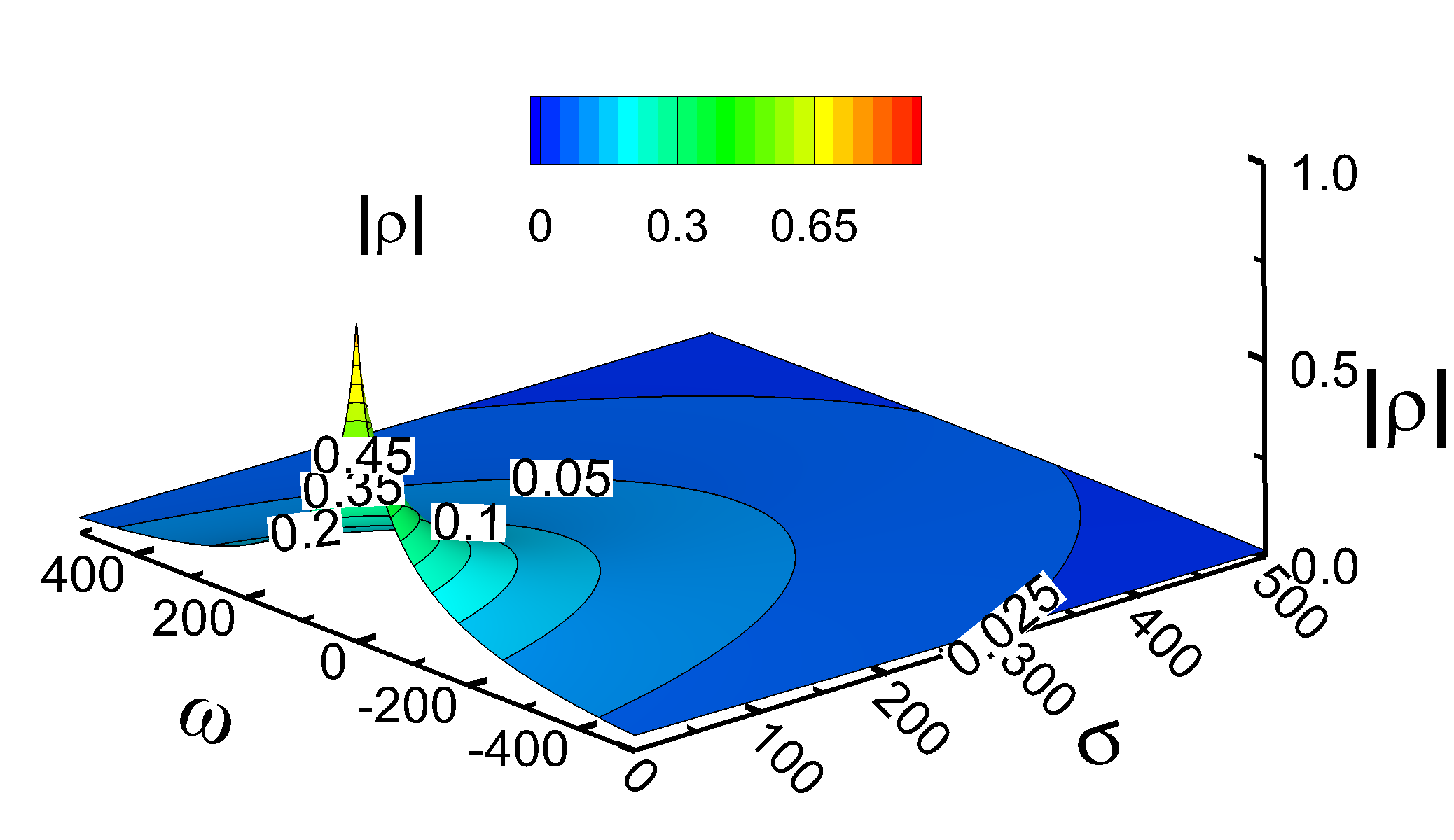}} \ 
  \subfloat[][]
  {\label{fig3e}\includegraphics[width=0.32\textwidth]{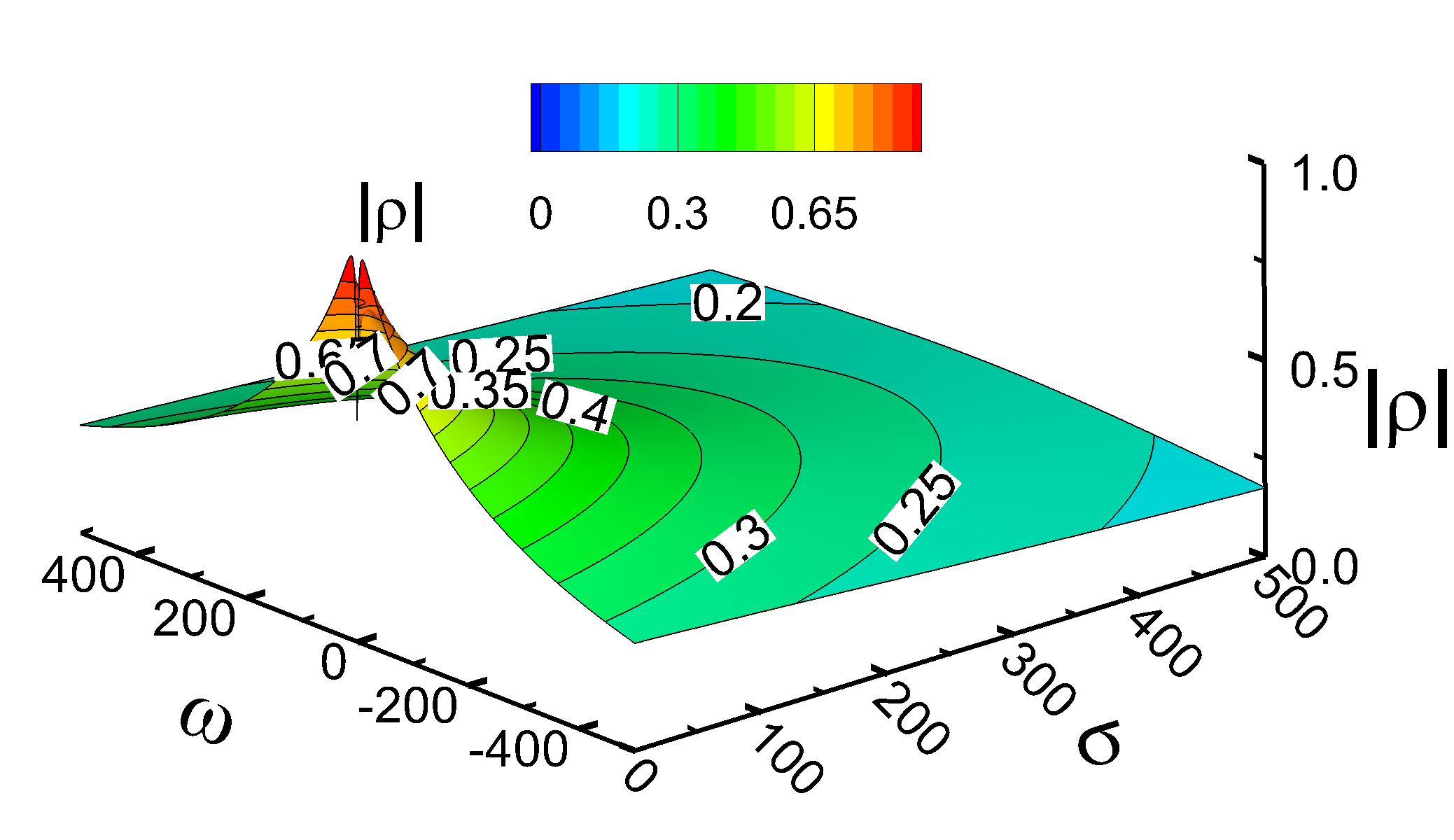}} \ 
  \subfloat[][]
  {\label{fig3f}\includegraphics[width=0.32\textwidth]{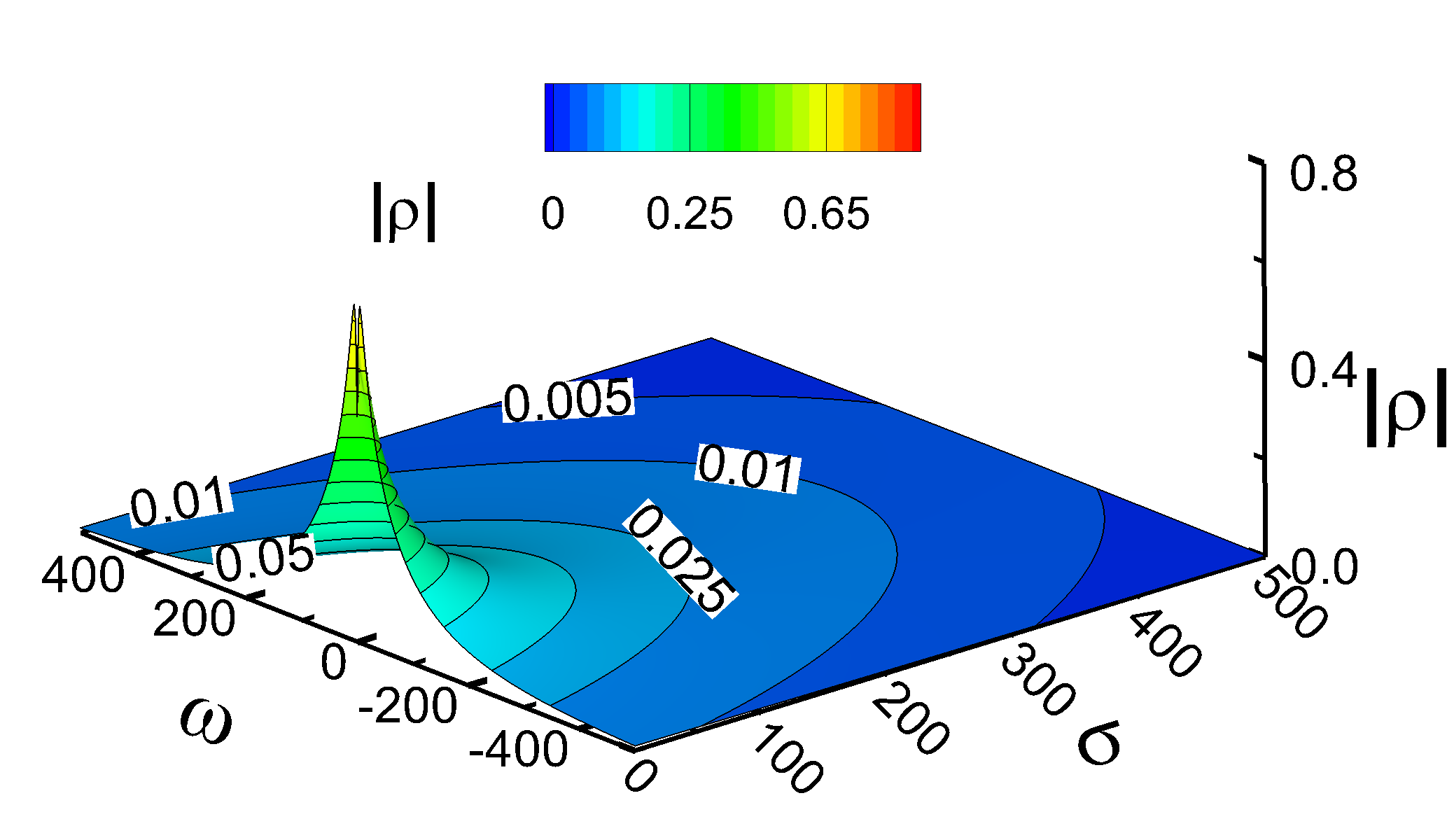}} \\ 
  \vspace{-8pt}
  \caption{\small Distributions of $|\rho |$ associated with different transmission conditions. 
  a) Condition (\ref {numerical_Dirichlet_cond}), or (\ref {numerical_Robin2_cond}), Case 1. 
  b) Condition (\ref {numerical_comb2_cond}), Case 1.  
  c) Condition (\ref {numerical_comb3_cond}), Case 1.  
  d) Condition (\ref {numerical_Dirichlet_cond}), or  (\ref {numerical_Robin2_cond}), Case 2.
  e) Condition (\ref {numerical_comb2_cond}), Case 2.
  f) Condition (\ref {numerical_comb3_cond}), Case 2.}
  \label{fig3}
  \vspace{-0pt}
\end{figure}

\ul {Numerical examples} Consider computation of two cases as shown in Table \ref {case12}. With regards to behaviors of solutions, Case 1 is advection-dominant, while Case 2 is diffusion-dominant. The contract factors are plotted in Fig. \ref {fig3} for the two cases associated with all above transmission conditions, except that of the Robin condition (\ref {numerical_Robin1_cond}) with two overlapping nodes, for which $\rho =1$. It is seen in Figs. \ref {fig3a} and \ref {fig3d} that, in the both cases, the factor has a maximum at $(\omega ,\sigma)=(0,0)$, when they are associated with the Dirichlet condition (\ref {numerical_Dirichlet_cond}) and the Robin condition with three overlapping nodes (\ref {numerical_Robin2_cond}). If the combined conditions (\ref {numerical_comb2_cond}) and (\ref {numerical_comb3_cond}) are adopted, there are one minimum and two maximums around $(\omega ,\sigma)=(0,0)$. From (\ref {rho_comb2}) and (\ref {rho_comb3}), it is seen that the difference of their contraction factors is that of the latter has a factor $r _-/r_+$, which is less then 1 in magnitude. Therefore, it is expected that distributions for the contraction factors with the two conditions have a similar shape, but the latter is smaller in magnitude. This is indeed the case, as seen in the comparison of Figs. \ref {fig3b} and \ref {fig3c}, and Figs. \ref {fig3e} and \ref {fig3f}. Note that, as indicated previously, $\rho $ is analytical, so its maximums/minimums only happen at $\sigma =0$. This is the case in all of these examples, as seen in Fig. \ref {fig3}. 

\begin{table}[!ht]
  \centering
  \caption{\small Cases for numerical experiments. $\Delta x=0.05$.} 
  \label{case12}
  \begin{tabular}{|c|c|c|}
\hline
       Case  & $\mu _1$, $\theta _1$, $\gamma _1$ & $\mu _2$, $\theta _2$, $\gamma _2$ \\ \hline
        1    & 1.5, \ 0.1, \ 0                 &  1, \ 0.1, \ 0   \\ \hline       
        2    & 0.2, \ 0.4, \ 1                 &  0.4, \ 0.2, \ 2 \\ \hline   
  \end{tabular}
\end{table}

\begin{figure}[!ht] 
  \vspace{-0.2 cm}
  \centering
  \subfloat[][]  
  {\label{fig4a}\includegraphics[width=0.32\textwidth]{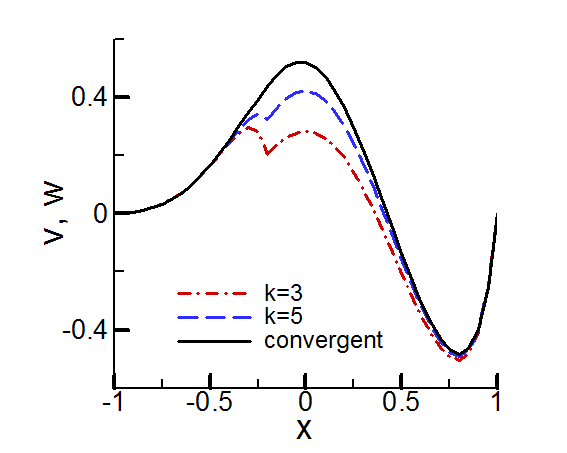}}  
  \subfloat[][]
  {\label{fig4b}\includegraphics[width=0.32\textwidth]{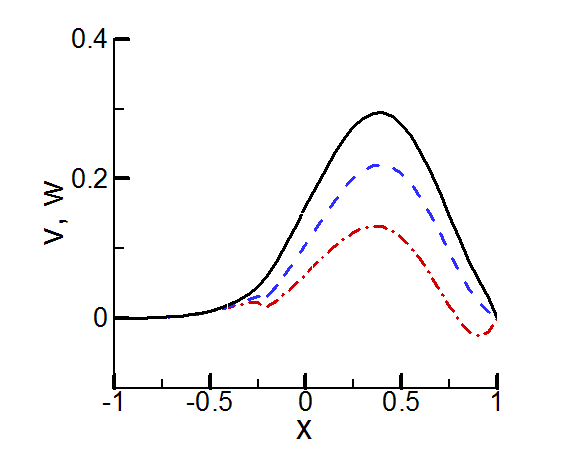}} 
  \subfloat[][]
  {\label{fig4c}\includegraphics[width=0.32\textwidth]{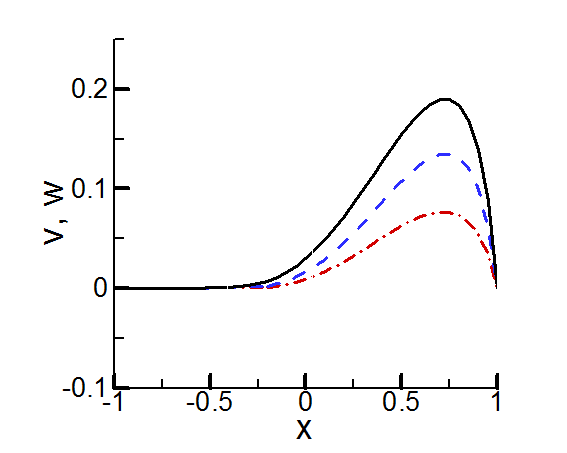}} \\ 
  \subfloat[][]
  {\label{fig4d}\includegraphics[width=0.32\textwidth]{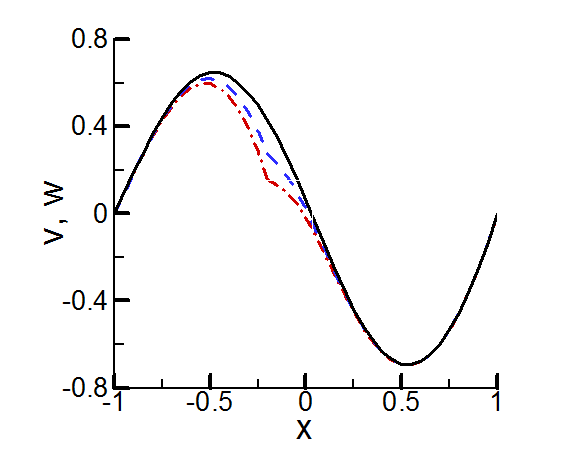}} 
  \subfloat[][]
  {\label{fig4e}\includegraphics[width=0.32\textwidth]{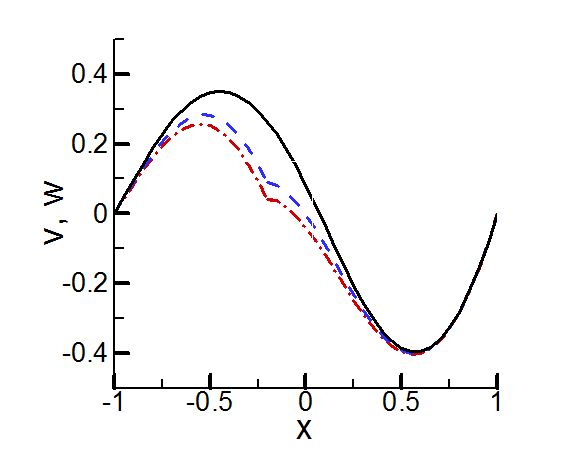}} 
  \subfloat[][]
  {\label{fig4f}\includegraphics[width=0.32\textwidth]{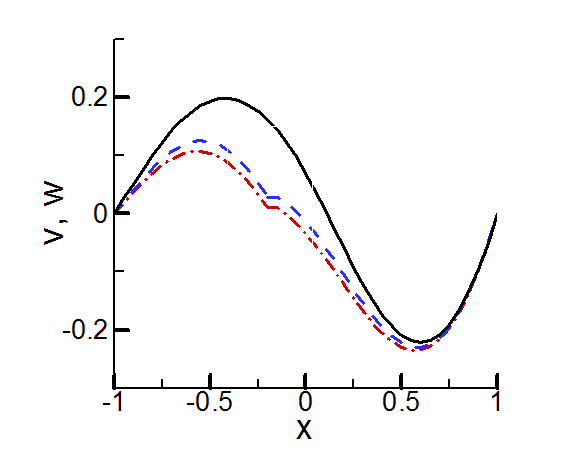}} 
  \vspace{-8pt}
  \caption{\small Solutions at different Schwarz iterations associated with the Dirichlet condition, (\ref {numerical_Dirichlet_cond}). a) Case 1, t=0.5. b) Case 1, t=1. c) Case 1, t=1.5. d) Case 2, t=0.15. e) Case 2, t=0.3. f) Case 2, t=0.45.}
  \label{fig4}
  \vspace{-0pt}
\end{figure}

Let us compute the coupled equations, Eq. (\ref  {dis_ivp1}), for the two cases using the Dirichlet transmission condition, i.e., Eq. (\ref {numerical_Dirichlet_cond}), on $x \in [-1,1]$, $t \in [0,5]$. The initial condition is $f(x)=-sin (\pi x)$, and boundary conditions $v=0$ at $x=-1$, $w=0$ at $x=1$ are adopted. The left subdomain is $-1\le x\le -0.15$, and the right subdomain is $-0.2\le x\le 1$. The left and right subdomain have an interface at $x_1=-0.15$, $x_2=-0.2$, respectively, and their solutions there are $v_0^k$ and $w_0^k$, respectively. For the computation, the backward Euler method is used for the discretization in time, with grid spacing $\Delta x=0.05$, and time step $\Delta t=0.05$. To make it simple, the initial values for the iteration at the interfaces are set as zero, that is $v_0^0=0$, $w_0^0=0$, $t\in [0,5]$. The solutions at three time instants computed at different iteration and convergence are plotted in Fig. \ref {fig4}. Here and hereafter, convergent solutions are the iterated values at fully convergence. The figure shows the convergence processes of the iterated solutions.

\section {Optimized transmission algorithm}

\ul {Optimized algorithm for Dirichlet condition} In order to expedite the convergence of Schwarz iteration between the two subdomains, we consider optimized transmission algorithms. For the Dirichlet condition (\ref {numerical_Dirichlet_cond}), an optimized algorithm is   
\begin{equation}
\label{opt_numerical_Dirichlet_cond}
\begin{array}{l}
(v_0^{k+1}-v_{-1}^{k+1})+\alpha  v_0^{k+1}=(w_1^k-w_0^k)+\alpha  w_1^k \\ 
(w_1^{k+1}-w_0^{k+1})+\beta w_0^{k+1}=(v_0^k-v_{-1}^k)+\beta v_{-1}^k 
\end{array}
\end{equation}

\noindent Such optimal algorithm was first proposed for waveform relaxation of circuit problems by Gander et al., and it was effective in speeding up the convergence  \cite {Gander:2004}. As they indicated, also it is easy to verify, above condition recovers to algorithm (\ref {numerical_Dirichlet_cond}) when the iteration converges as long as 
\begin{equation}
\label{int_cond_recovery}
(\alpha +1)(\beta -1)+1\ne 0. 
\end{equation}

\noindent Algorithm (\ref {opt_numerical_Dirichlet_cond}) leads to 
\begin{equation*}
\begin{array}{l}
v_0^{k+1}=\dfrac {1}{1+\alpha  }v_{-1}^{k+1}-\dfrac {1}{1+\alpha }w_0^k+w_1^k \\
w_0^{k+1}=\dfrac {1}{1-\beta }w_1^{k+1}+v_{-1}^k-\dfrac {1}{1-\beta }v_0^k
\end{array}
\end{equation*}

\noindent Then, by the same steps in the above section, one derives that 
\begin{equation*}
\begin{array}{l}
\rho =\dfrac {c_1{r_-}^{-1}-(1+\alpha )c_1}{c_1+(1+\alpha )(a_1{r_+}^{-1}+b_1-s)} \cdot 
      \dfrac {a_2r_+-(1-\beta )a_2}{a_2+(1-\beta )(c_2r_-+b_2-s)}
\end{array}
\end{equation*}

\noindent which is further derived as 
\begin{equation}
\label{opt_rho_Dirichlet}
\begin{array}{l}
\rho =\dfrac {1+\alpha  -{r_-}^{-1}}{1+\alpha  -{r_+}^{-1}} \cdot
      \dfrac {1-\beta -r_+}{1-\beta -r_-} \cdot 
      \dfrac {r_-}{r_+}
\end{array}
\end{equation}

\vspace {0.3 cm}

{\textsc{Theorem 4.1}} Let (\ref {dis_ivp1}) be computed via the optimized algorithm (\ref {opt_numerical_Dirichlet_cond}). Then,  i) $\rho $ in (\ref {opt_rho_Dirichlet}) is an analytical function under condition  
\begin{equation}
\label{cond_alpha_beta}
\alpha >0, \ \beta <0  
\end{equation}

\noindent ii) A set of optimal values for $(\alpha ,\beta)$ is 
\begin{equation}
\label{opt_alpha_beta_Dirichlet}
\alpha ^\ast =-1+{r_-}^{-1}, \ \ \beta ^\ast =1-r_+
\end{equation}

\noindent which leads to $\rho =0$ in (\ref {opt_rho_Dirichlet}).

\vspace {0.3 cm}

{\textsc{Proof}}  First, since roots $r_{\pm}$ are analytical, it is suffices to show that $\rho $ is analytical if we can prove that the denominator of $\rho $ in (\ref {opt_rho_Dirichlet}) is not zero. In view of (\ref {r_-+}), $r_+\ne 0$. Additionally, because of (\ref {r_-+}) and (\ref {cond_alpha_beta}),  $1+\alpha -r_+^{-1}$, ${1-\beta -r_-}$ $\ne 0$. Therefore, the denominator is not zero.

Second,  in view of (\ref {r_-+}),  $r_-$ and $r_+$ are different in magnitude. As a result, while condition (\ref {opt_alpha_beta_Dirichlet}) does not make that the denominator of (\ref {opt_rho_Dirichlet}) to become zero, it leads to that the numerator of (\ref {opt_rho_Dirichlet}) becomes zero. Therefore, under condition (\ref {opt_alpha_beta_Dirichlet}), $\rho =0$. This completes the proof $ \sharp $. 

\vspace {0.3 cm}

{\textsc{Remark 4.1}} The optimal values $(\alpha ^\ast,\beta ^\ast)$ lead to the best possible values for contraction factor, i.e., $\rho =0$, and this will make the computation to achieve convergence within two times of iteration, given that $\rho =\hat {v}_i^2/\hat {v}_i^0, \ \hat {w}_i^2/\hat {w}_i^0$, see (\ref {Laplace_iter_conv}) (here superscripts $0$ and $2$ mean iteration induces). However, the optimal values $(\alpha ^\ast,\beta ^\ast)$ are difficult to implement in practice, because they are complex numbers and non-constant in the time domain. Furthermore, since $\rho $ in (\ref {opt_rho_Dirichlet}) is analytical, it achieves maximums and maximums only on the boundary of the domain, that is, on $s =\omega i$.     

\vspace {0.3 cm}

{\textsc{Remark 4.2}} The conclusions in above theorem have been obtained in a study on the RC-type circuits, in which $a_1=a_2=c_1=c_2$, $b_1=b_2$ \cite {Gander:2004}, and they may be considered as an extension of the previous work.  

\vspace {0.3 cm}

\ul {Optimized algorithm for Robin condition} Following the idea of algorithm (\ref {opt_numerical_Dirichlet_cond}), an optimized algorithm for flux condition (\ref {numerical_Robin2_cond}) is constructed as 
\begin{equation}
\label{opt_numerical_Robin2_cond}
\begin{array}{l}
(\hat f_{-1/2}^{k+1}-\hat f_{-3/2}^{k+1})+\alpha  \hat  f_{-1/2}^{k+1}=(\hat g_{3/2}^k-\hat g_{1/2}^k)+\alpha  \hat g_{3/2}^k \\
(\hat g_{3/2}^{k+1}-\hat g_{1/2}^{k+1})+\beta \hat g_{1/2}^{k+1}=(\hat f_{-1/2}^k-\hat f_{-3/2}^k)+\beta \hat f_{-3/2}^k
\end{array}
\end{equation}

\noindent It can be verified that the algorithm recovers to (\ref {numerical_Robin2_cond}) at convergence, as long as (\ref {int_cond_recovery}) holds. From (\ref {opt_numerical_Robin2_cond}), one has 
\begin{equation*}
\begin{array}{l}
v_0^{k+1}=-\dfrac {a_1}{(1+\alpha )c_1}v_{-2}^{k+1}+\dfrac {(1+\alpha )a_1+c_1}{(1+\alpha )c_1}v_{-1}^{k+1}
          +\dfrac {a_2}{(1+\alpha)c_1}w_0^k-\dfrac {(1+\alpha)a_2+c_2}{(1+\alpha )c_1}w_1^k+\dfrac {c_2}{c_1}w_2^k \\
w_0^{k+1}=\dfrac {a_2+(1-\beta )c_2}{(1-\beta )a_2}w_1^{k+1}-\dfrac {c_2}{(1-\beta )a_2}w_2^{k+1}
         +\dfrac {a_1}{a_2}v_{-2}^k-\dfrac {a_1+(1-\beta )c_1}{(1-\beta )a_2}v_{-1}^k+\dfrac {c_1}{(1-\beta )a_2}v_0^k 
\end{array}
\end{equation*}

\noindent In this situation, the contraction factor is derived as   
\begin{equation*}
\label{opt_rho_Robin2}
\begin{array}{lll}
\rho
&=&\dfrac {c_2-a_2{r_-}^{-1}+(1+\alpha )(a_2-c_2r_-)}{c_1-a_1{r_+}^{-1}+(1+\alpha )(a_1-c_1r_+)} \\
&&\cdot \dfrac {a_1r_+-c_1{r_+}^2+(1-\beta)(c_1r_+-a_1)}{a_2r_--c_2{r_-}^2+(1-\beta )(c_2r_--a_2)}\cdot \dfrac {r_-}{r_+}  \\ 
&=&\dfrac {(c_2r_--a_2){r_-}^{-1}+(1+\alpha )(a_2-c_2r_-)}{(c_1r_+-a_1){r_+}^{-1}+(1+\alpha )(a_1-c_1r_+)} \\
&&\cdot \dfrac {(a_1-c_1r_+)r_++(1-\beta)(c_1r_+-a_1)}{(a_2-c_2r_-)r_-+(1-\beta )(c_2r_--a_2)}\cdot \dfrac {r_-}{r_+}  \\ 
&=&\dfrac {(1+\alpha -{r_-}^{-1})(a_2-c_2r_-)}{(1+\alpha -{r_+}^{-1})(a_1-c_1r_+)} \\
&&\cdot \dfrac {(-1+\beta +r_+)(a_1-c_1r_+)}{(-1+\beta +r_-)(a_2-c_2r_-)}\cdot \dfrac {r_-}{r_+}  \\ 
&=&\dfrac {1+\alpha  -{r_-}^{-1}}{1+\alpha  -{r_+}^{-1}} \cdot
      \dfrac {1-\beta -r_+}{1-\beta -r_-} \cdot 
      \dfrac {r_-}{r_+} 
\end{array}
\end{equation*}

\noindent This indicates that the contraction factor is same to that in (\ref {opt_rho_Dirichlet}).

\vspace {0.3 cm}

{\textsc{Theorem 4.2}} Let (\ref {dis_ivp1}) be computed via the optimized algorithm for Robin condition (\ref {opt_numerical_Robin2_cond}). Then, i) the contraction factor is the same to that via the optimized algorithm (\ref {opt_numerical_Dirichlet_cond}), i.e., (\ref {opt_rho_Dirichlet}). ii) The contraction factor is an analytical function under condition (\ref {cond_alpha_beta}). iii) A set of optimal values for $(\alpha ,\beta)$ is (\ref {opt_alpha_beta_Dirichlet}). 

\vspace {0.3 cm}

{\textsc{Proof}} The proof is the same to that of Theorem 4.1, since their contraction factors are the same.

\vspace {0.3 cm}

\ul {Numerical experiments} 
Now, let us look for optimal values of the parameters, ($\alpha ^\ast $,$\beta ^\ast$).  Particularly, to reduce the magnitude of the contraction factor in its whole range along $\sigma =0$, the optimal values are so chosen that the following is satisfied: 
\begin{equation}
\label{optimal1}
(\alpha ^\ast ,\beta ^\ast)=\lbrace \left ( \alpha ,\beta  \right ) \big\rvert |\rho |=min_{\alpha ,\beta } (max _{Re(s)>0} |\rho (\alpha ,\beta ,s,a,b,c)|) \rbrace
\end{equation}

\noindent  which requires that the contraction factor remains small in magnitude in all range of $s$ so convergence speeds up. In view that the contraction factor is a function of $s$, in general it is difficult to solve above for the optimal values ($\alpha ^\ast $,$\beta ^\ast$). Actually, we just need to find these values along $s=\omega i$, where the minimums and maximums of $\rho $ reside. As a practical approach, we follow \cite {Gander:2004} by letting $\alpha ^\ast =-\beta ^\ast $ and solving $\alpha ^\ast $ from 
\begin{equation}
\label{optimal2}
\rho (\alpha ^\ast ,-\alpha ^\ast ,i\omega _{min} ,a,b,c)=\rho (\alpha ^\ast ,-\alpha ^\ast ,i\omega _{max},a,b,c)
\end{equation} 

\noindent Here $\omega _{min}$ and $\omega _{max} $ may be determined as $\pi /T$ and $\pi /\Delta t$, with $T $ and $\Delta t$ being the time length of problem (\ref {ivp1}) and the small time scale to be resolved, respectively. 

\begin{figure}[!ht] 
  \vspace{-0.2 cm}
  \centering
  \subfloat[][]  
  {\label{fig5a}\includegraphics[width=0.32\textwidth]{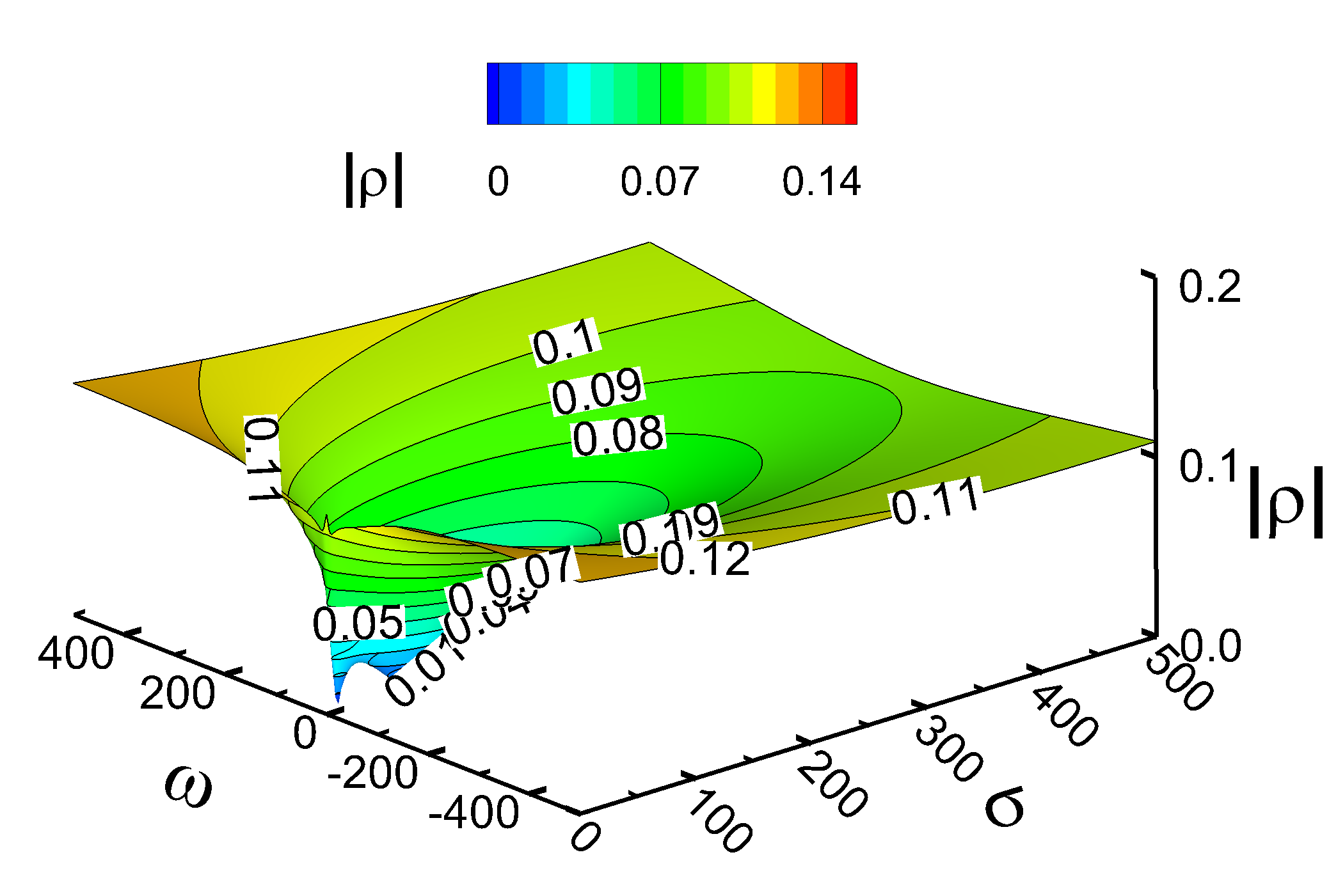}} \ \
  \subfloat[][]  
  {\label{fig5b}\includegraphics[width=0.32\textwidth]{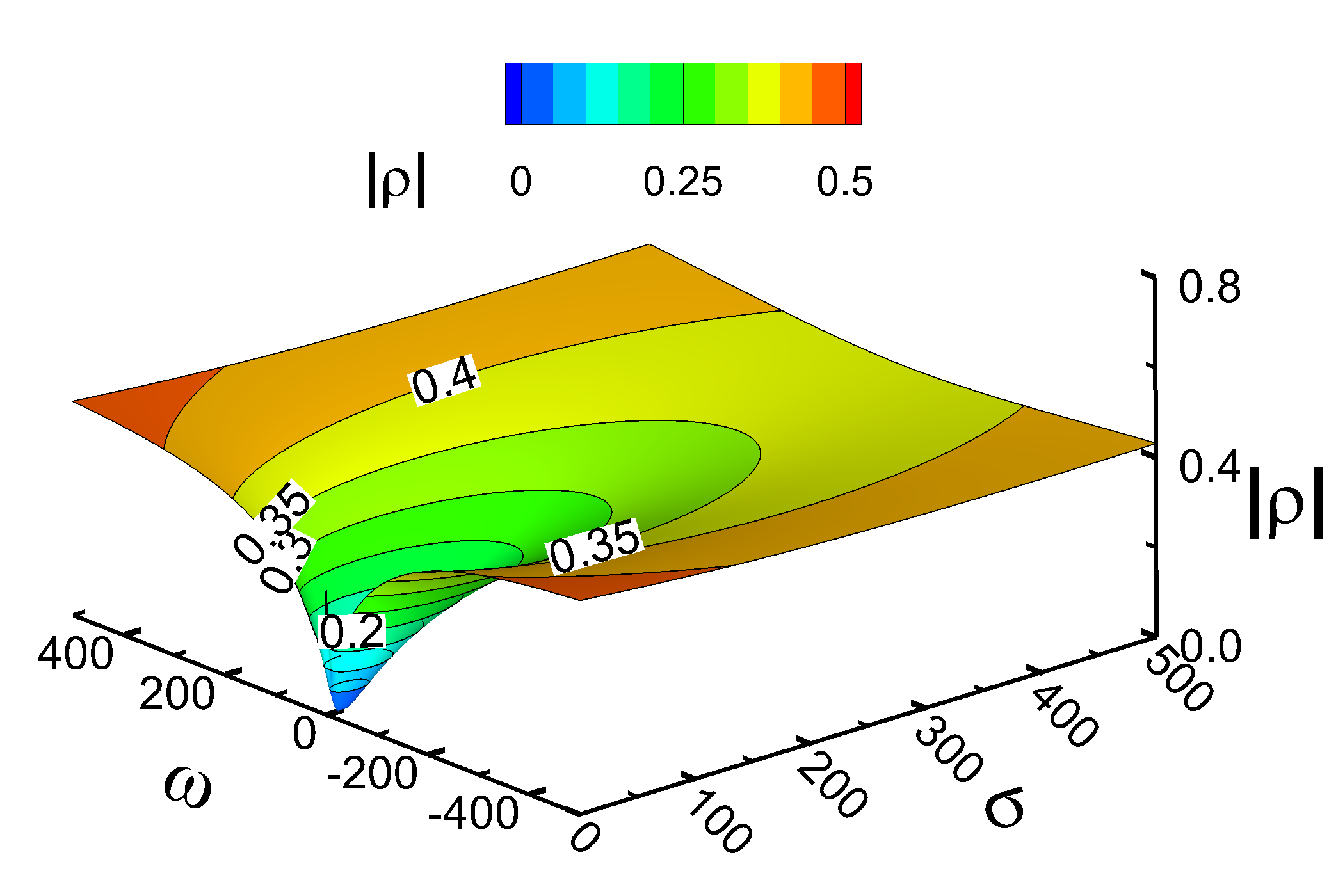}}  
  \vspace{-8pt}
  \caption{\small Distributions of $|\rho |$ (\ref {opt_rho_Dirichlet}). a) Case 1. $\alpha ^*=1.76848$. b) Case 2. $\alpha ^*=0.40361$.}
  \label{fig5}
  \vspace{-0pt}
\end{figure}

With the above procedure, the optimal value $\alpha ^*$ is estimated, and the corresponding $|\rho |$ for Case 1 and 2 is plotted in Fig. \ref {fig5}. By comparison of Figs. \ref {fig3a} and \ref {fig5a}, \ref {fig3d} and \ref {fig5b}, it is seen that $|\rho |$ of the optimized algorithm has changed substantially in terms of distribution from that of the original algorithm. Additionally, its overall values, including the minimums and maximums at the spikes and trough, become smaller in magnitude. 

\begin{figure}[!ht] 
  \vspace{-0.2 cm}
  \centering
  \subfloat[][]  
  {\label{fig6a}\includegraphics[width=0.32\textwidth]{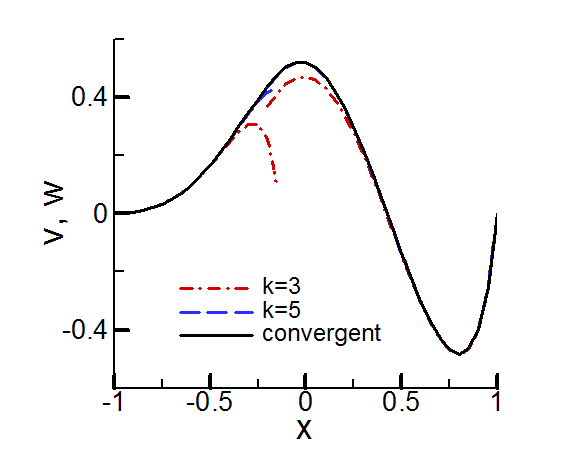}}  
  \subfloat[][]
  {\label{fig6b}\includegraphics[width=0.32\textwidth]{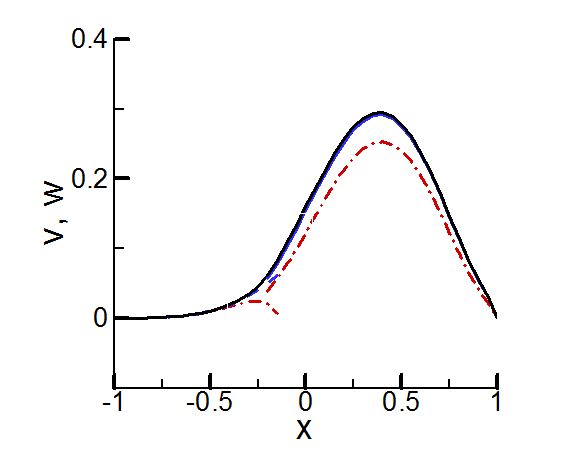}} 
  \subfloat[][]
  {\label{fig6c}\includegraphics[width=0.32\textwidth]{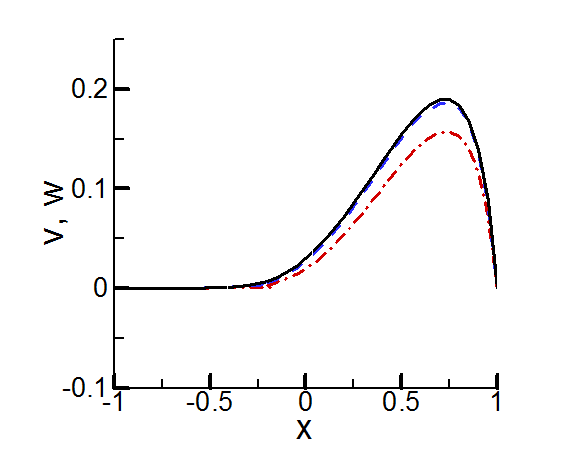}} \\ 
  \subfloat[][]
  {\label{fig6d}\includegraphics[width=0.32\textwidth]{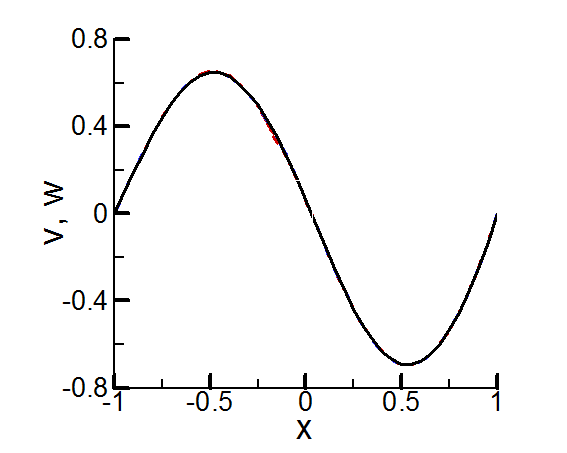}} 
  \subfloat[][]
  {\label{fig6e}\includegraphics[width=0.32\textwidth]{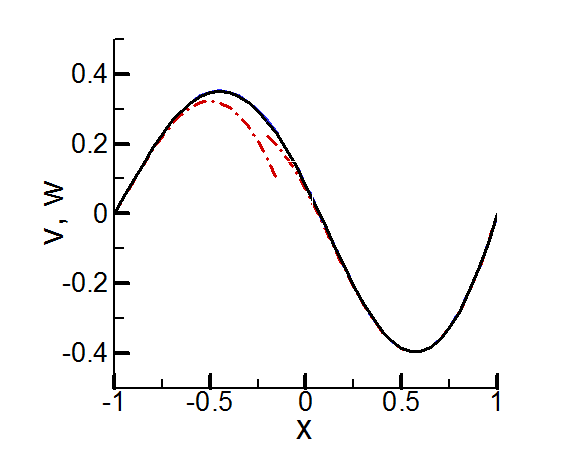}} 
  \subfloat[][]
  {\label{fig6f}\includegraphics[width=0.32\textwidth]{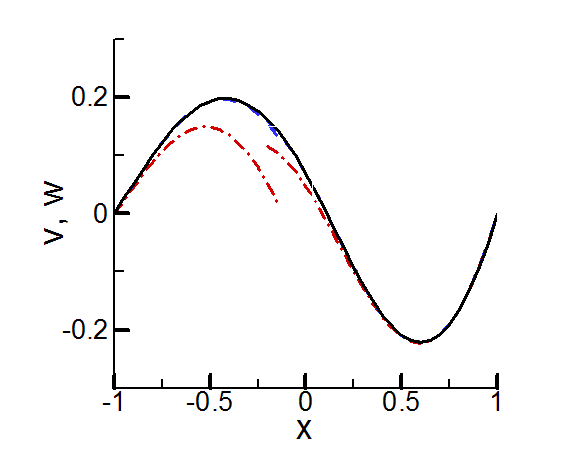}} 
  \vspace{-8pt}
  \caption{\small Solutions at different Schwarz iterations associated with the optimized condition, (\ref {opt_numerical_Dirichlet_cond}). a) Case 1, t=0.5. b) Case 1, t=1. c) Case 1, t=1.5. d) Case 2, t=0.15. e) Case 2, t=0.3. f) Case 2, t=0.45.}
  \label{fig6}
  \vspace{-0pt}
\end{figure}

\begin{figure}[!ht] 
  \vspace{-0.2 cm}
  \centering
  \subfloat[][]  
  {\label{fig7a}\includegraphics[width=0.32\textwidth]{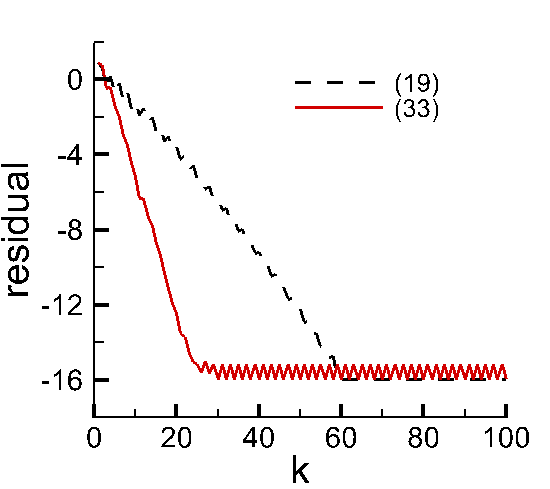}} \
  \subfloat[][]  
  {\label{fig7b}\includegraphics[width=0.32\textwidth]{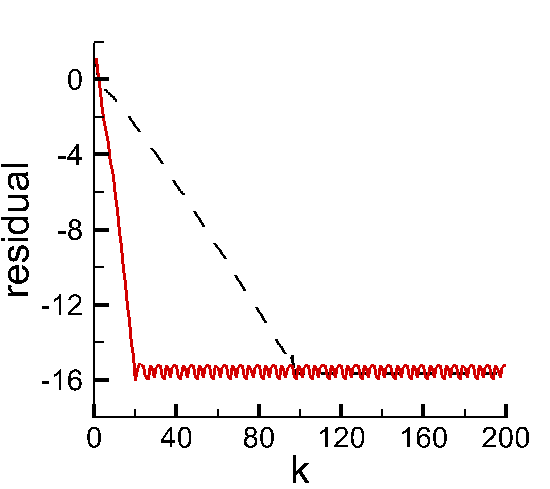}} 
  \vspace{-8pt}
  \caption{\small Convergence of the computation with Conditions (\ref {numerical_Dirichlet_cond}) and (\ref {opt_numerical_Dirichlet_cond}). a) Case 1, $t=0.5$. b) Case 2, $t=0.15$.}
  \label{fig7}
  \vspace{0pt}
\end{figure}

Let us compute the problem with the optimized transmission condition, i.e., Eqs. (\ref  {dis_ivp1}) and (\ref {opt_numerical_Dirichlet_cond}), with same parameters as before. Again, the computed solutions at different Schwarz iteration and convergence are plotted in Fig. \ref {fig6}. In comparison to Fig. \ref {fig4}, it is seen that, if the optimized transmission algorithm is adopted, the iterated solutions reach the convergent solutions much faster in terms of iteration number. Moreover, it is noticed that the iterated solutions may not be continuous at the subdomain interfaces before the full convergence. To further compare the convergence speed, samples of the iteration residuals in terms of 1-norm are plotted in Fig. \ref {fig7}. The figure shows that indeed the optimized algorithm speeds up the convergence of the computation.  

\section {Extension to nonlinear equations}

Now, let us consider computation of coupled nonlinear equations via Schwarz iteration as follows:   
\begin{equation}
\label{ivp1_nonlinear}
\begin{array}{ll}
\left\{\begin{array}{ll}
v_t^{k+1} + \mu (v^k)v_x^{k+1} = \theta _1v_{xx}^{k+1}-\gamma \times (v^k)^\kappa v^{k+1},  \ & t \in (0,T] \\
v^{k+1}=f(x),  \ & t=0 \\
v^{k+1}=p_1(w^k), \ & x=x_2 
\end{array} \right.  \\   [18pt]
\left\{\begin{array}{ll}
w_t^{k+1} + \mu (w^k)w_x^{k+1} = \theta _2w_{xx}^{k+1}-\gamma \times (w^k)^\kappa w^{k+1}, \ & t \in (0,T]  \\
w^{k+1}=g(x), \ & t=0 \\ 
w^{k+1}=p_2(v^k), \ & x=x_1 
\end{array} \right. 
\end{array}
\end{equation}
\noindent Here, $\mu (\cdot )$ is a function of the solutions. $\kappa =0$, $1$, which corresponds to the first- and second-order reaction, respectively, and they are common scenarios in practice, e.g., in water flows \cite {Ji:2008}. The nonlinearity of the above equations comes from the advection and source terms. As an extension to that of the linear equations, i.e., (\ref {dis_ivp1}), the above problem is discretized as follows 
\begin{equation}
\label{dis_ivp1_nonlinear}
\begin{array}{ll}
\left\{\begin{array}{ll}
{v_t}_i^{k+1}=a_{1,i}^kv_{i-1}^{k+1}+b_{1,i}^kv_i^{k+1}+c_{1,i}^kv_{i+1}^{k+1},  \ \ i\le -1 \\
v_0^{k+1}=\tilde {p}_1(w_0^k,w_1^k, ...)  \\
\end{array}\right . \\ [10pt]
\left\{\begin{array}{ll}
{w_t}_i^{k+1}=a_{2,i}^kw_{i-1}^{k+1}+b_{2,i}^kw_i^{k+1}+c_{2,i}^kw_{i+1}^{k+1},  \ \ i\ge 1 \\
w_0^{k+1}=\tilde {p}_2(v_0^k,v_{-1}^k, ...)\\ 
\end{array}\right . 
\end{array}
\end{equation}
\noindent in which   
\begin{equation*}
\label{abc_noinlinear}
%\left \{
\begin{array}{lll}
a_{1,i}^k=\dfrac {\mu (v^k_i)}{2\Delta x}+\dfrac {\theta }{{\Delta x}^2}, \ \
& b_{1,i}^k=-\dfrac {2\theta }{{\Delta x}^2}-\gamma \times (v^k_i)^{\kappa _1}, \ \ 
& c_{1,i}^k=-\dfrac {\mu (v^k_i)}{2\Delta x}+\dfrac {\theta }{{\Delta x}^2} \\
a_{2,i}^k=\dfrac {\mu (w^k_i)}{2\Delta x}+\dfrac {\theta }{{\Delta x}^2}, \ \
& b_{2,i}^k=-\dfrac {2\theta }{{\Delta x}^2}-\gamma \times (w^k_i)^{\kappa _2}, \ \ 
& c_{2,i}^k=-\dfrac {\mu (w^k_i)}{2\Delta x}+\dfrac {\theta }{{\Delta x}^2}
\end{array} 
%\right.
\end{equation*}
\noindent 

In the above computational problem, because coefficients $a_{l,i}^k$, $b_{l,i}^k$, and $c_{l,i}^k$ ($l=1,2$) are not constants but solution dependent, we are unable to directly derive contraction factors as for the linear equations. Before discussion on contraction factors for the problem, we consider a property of contraction factors for the linear equations. Particularly, from (\ref {int_Sch1_solu}), it is derived that 
\begin{equation*}
\begin{array}{l}
|\hat {v}_i^{k+1}|=|B^{k+1}{r_+}^i|=|B^{k+1}{r_+}^{i-1}{r_+}|=|\hat {v}_{i-1}^{k+1}||r_+|
= \ ... \ =|\hat {v}_{i-q}^{k+1}||r_+|^q,  \ \  i=0, -1, -2, ... \\
|\hat {w}_{i+1}^{k+1}|=|A^{k+1}{r_-}^{i+1}|=|A^{k+1}{r_-}^i{r_-}|=|\hat {w}_i^{k+1}||r_-|
= \ ... \ =|\hat {w}_{i-q}^{k+1}||r_-|^{q+1},  \ \  i=0, 1,  2, ...
\end{array}
\end{equation*}

\noindent where $q $ is a positive integer. In view of (\ref {r_-+}), the above leads to the following property.

\vspace {0.3 cm}

{\textsc{proposition 5.1}}  The solution of (\ref {trans_int_Sch}) exhibits the following property,   
\begin{equation}
\label{int_Sch1_solu_prop}
\begin{array}{l}
|\hat {v}_i^{k+1}|<|\hat {v}_{i+1}^{k+1}|\cdots <|\hat {v}_0^{k+1}|,  \ \  i=0, -1, -2, ... \\
|\hat {w}_0^{k+1}|>\cdots |\hat {w}_i^{k+1}|>|\hat {w}_{i+1}^{k+1}|,  \ \  i=0, 1,  2, ...
\end{array}
\end{equation}

\vspace {0.3 cm}

Proposition 5.1 indicates that the magnitude of  $ |\hat {v}_i^{k+1}| $ and $ |\hat {w}_i^{k+1}| $ decade as they are located away from the interface locations (at $i=0$), or, as $|i|$ increases. As seen in the above, such decaying is fast in view that $|r_-|<1$, $|r_+|>1$. These imply that the the values of $ |\hat {v}_0^{k+1}| $ and $ |\hat {w}_0^{k+1}| $ at the interfaces affect the values of $\hat {v}_i^{k+1} $ and $\hat {w}_i^{k+1} $ only in a small neighborhood near the interfaces, or, when $i$ is relatively small. In general, because of the nonlinearity, coefficients $a_{l,i}^k$, $b_{l,i}^k$, and $c_{l,i}^k$ in (\ref {dis_ivp1_nonlinear}) depend on solutions,  which in general changes with $x$, or, $i$. However, since effects of $ |\hat {v}_0^{k+1}| $ and $ |\hat {w}_0^{k+1}| $ at the interfaces are restricted to a small neighborhood of them, therefore, it is anticipated that the analysis and derivation for the contraction factors and optimized transmission conditions derived from the linear situations may be valid in certain degree. 

\begin{table}[!ht]
  \centering
  \caption{\small Cases for numerical experiments for nonlinear equations. } 
  \label{case34567}
  \begin{tabular}{|c|c|c|}
\hline
Case  & $\mu (v)$, $\theta _1$, $\gamma $, $\kappa _1$  & $\mu (w)$, $\theta _2$, $\gamma $, $\kappa _2$   \\ \hline
 3    & $v$,           \ 0.4, \ 1, \ 0     &  $w$,         \ 0.2,   \ 2, \ 0    \\ \hline 
 4    & $v^2$,         \ 0.4, \ 1, \ 0     &  $w^2$,       \ 0.2,   \ 2, \ 0    \\ \hline 
 5    & $sin(\pi v)$,  \ 0.4, \ 1, \ 0     &  $sin(\pi w)$,\ 0.2,   \ 2, \ 0    \\ \hline 
 6    & $v$,           \ 0.4, \ 1, \ 1     &  $w$,         \ 0.2,   \ 2, \ 1    \\ \hline 
 7    & $v$,           \ 0.04,\ 0, \ 0     &  $w$,         \ 0.02,  \ 0, \ 0    \\ \hline 
  \end{tabular}
\end{table}

With the above understanding, we extend the optimized transmission algorithm, i.e., (\ref {opt_numerical_Dirichlet_cond}), which is derived for Dirichlet condition (\ref {numerical_Dirichlet_cond}), to computation of (\ref {dis_ivp1_nonlinear}) by localization or linearization at the interfaces. Particularly, the roots are evaluated at the interfaces:
\begin{equation}
\label{roots3}
r_- =\frac {s-b_{2,0}^k - \sqrt {(s-b_{2,0}^k)^2-4a_{2,0}^kc_{2,0}^k}}{2c_{2,0}^k}, \quad  
r_+ =\frac {s-b_{1,0}^k + \sqrt {(s-b_{1,0}^k)^2-4a_{1,0}^kc_{1,0}^k}}{2c_{1,0}^k}
\end{equation}
\noindent In the experiments, similarly as in the linear scenarios, $f(x)=-sin (\pi x)$, $x \in [-1,1]$, $t \in [0,1]$. $u=0$ at $x=-1$, $v=0$ at $x=1$. $\Delta t=0.05$, $\Delta x=0.05$. $ v^k_0 $ and $w^k_0$ are located at  interfaces $x=-0.5$ and $x=-0.45$, respectively. Again, to make it simple, the initial values for the iteration at the interfaces are set as zero, that is $v_0^0=0$ $w_0^0=0$, $t\in [0,1]$.   

\begin{figure}[!ht] 
  \vspace{-0.2 cm}
  \centering
  \subfloat[][]  
  {\label{fig8a}\includegraphics[width=0.32\textwidth]{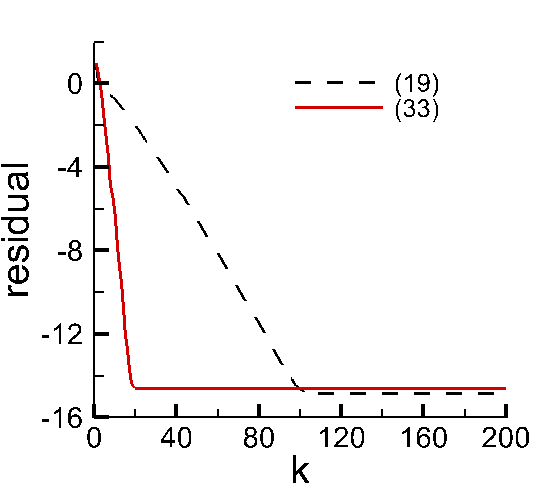}} \ 
  \subfloat[][]
  {\label{fig8b}\includegraphics[width=0.32\textwidth]{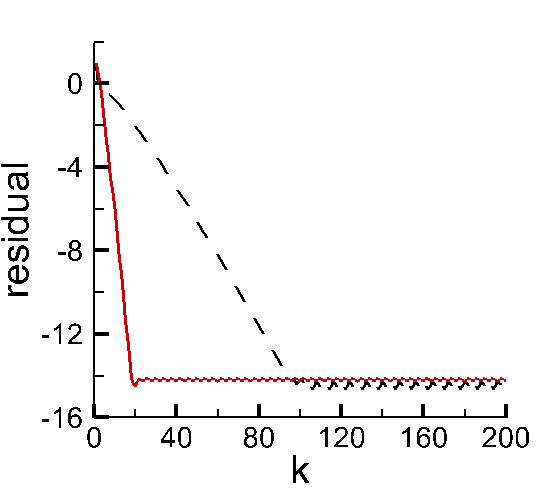}} \\ 
  \subfloat[][]
  {\label{fig8c}\includegraphics[width=0.32\textwidth]{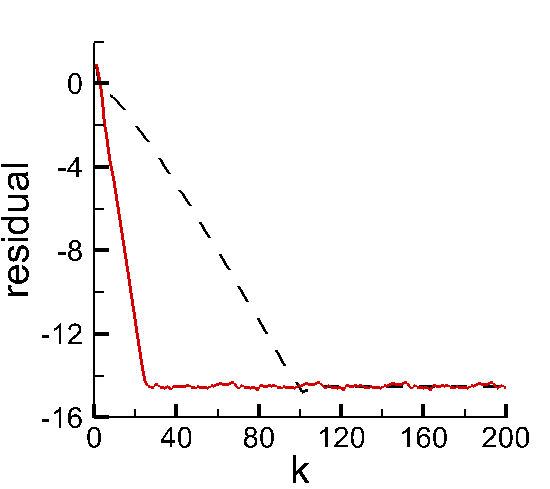}} \
  \subfloat[][]
  {\label{fig8d}\includegraphics[width=0.32\textwidth]{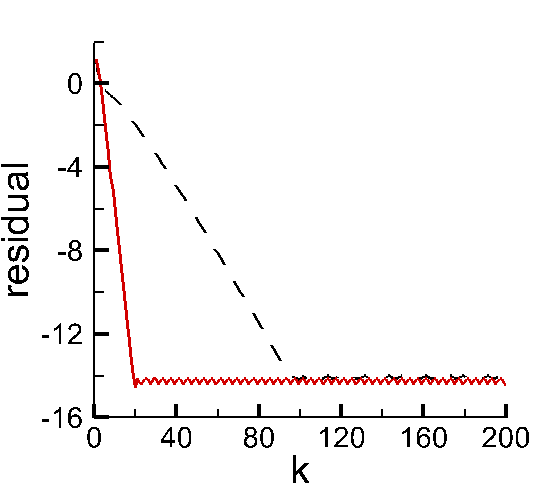}} \
  \subfloat[][]
  {\label{fig8e}\includegraphics[width=0.32\textwidth]{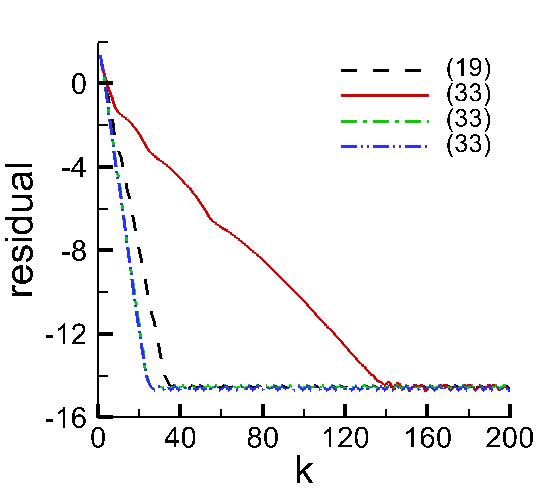}} \  
  \vspace{-8pt}
  \caption{\small Convergence of computation for nonlinear equations with Conditions (\ref {numerical_Dirichlet_cond}) and (\ref {opt_numerical_Dirichlet_cond}). $t=0.2$. a) Case 3, $\alpha ^*=0.48360$. b)  Case 4, $\alpha ^*=0.48436$. c) Case 5, $\alpha ^*=0.49372$. d) Case 6, $\alpha ^*=0.48292$. e) Case 7, $\alpha ^*_1=0.14600$, $\alpha ^*_2=2.05260$, $\alpha ^*_3=4.30670$.}
  \label{fig8}
  \vspace{0pt}
\end{figure}

Fig. \ref {fig8} presents experiments for five cases, which are detailed in Table \ref {case34567}, and it shows the convergence of the Schwarz iteration at an instant. In the computation, the method to search for the optimal values of $\alpha ^*$ is same as before. Interestingly, it is seen that in all of the first four cases (\ref {fig8a}, \ref {fig8b}, \ref {fig8c}, \ref {fig8d}), the optimized condition leads to much faster convergence in comparison to the Dirichlet condition, and the speedup of convergence is comparable to that in the linear cases in Fig. \ref {fig7}. 

However, the validity may fail as nonlinearity is involved, and Fig. \ref {fig8e} presents such an example. In this example, the method to search for optimal $\alpha $, particularly (\ref {optimal2}), leads to three values for $\alpha ^*$. The convergence associated with these values is plotted in the figure. In comparison to the Dirichlet condition, the optimized condition leads to little acceleration or even slowdown in convergence (two values for $\alpha ^*$ produce an almost same convergence) (Fig. \ref {fig8e}). 

\begin{figure}[!ht] 
  \vspace{-0.2 cm}
  \centering
  \subfloat[][]  
  {\label{fig9a}\includegraphics[width=0.32\textwidth]{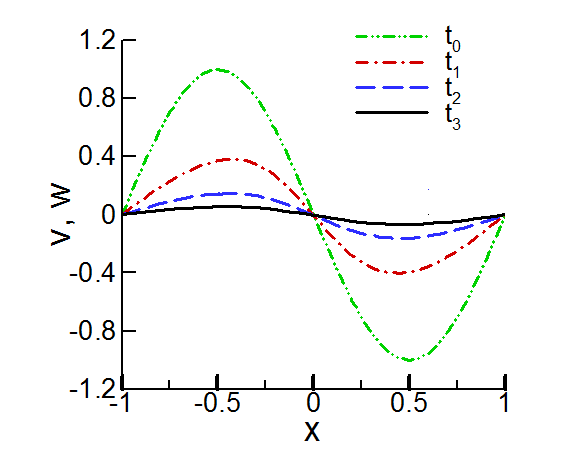}} \ 
  \subfloat[][]
  {\label{fig9b}\includegraphics[width=0.32\textwidth]{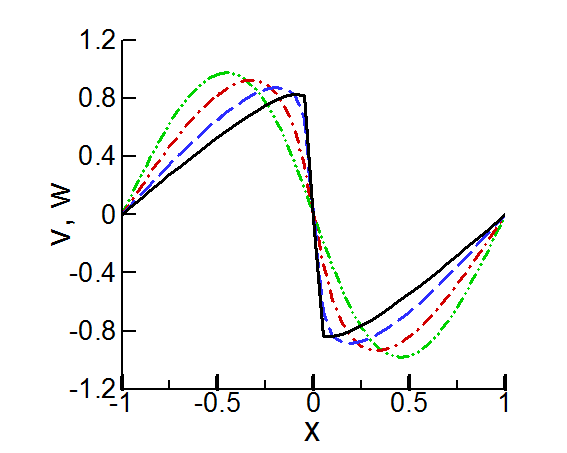}}
  \vspace{-8pt}
  \caption{\small Solutions for nonlinear equations at different Schwarz iterations and convergence. $t_0=0$, $t_1=0.2$, $t_2=0.4$, $t_3=0.6$.  a) Case3. b) Case 7.}
  \label{linear_nonlinear_solu}
  \vspace{0pt}
\end{figure}

The examples in Figs. \ref {fig8d} and \ref {fig8e} have a same type of equations but adopt different coefficients. As seen in Fig. \ref {linear_nonlinear_solu}, while the solution for Fig. (\ref {fig8d}) decays smoothly with time, the solution for Fig. (\ref {fig8e}) forms a discontinuity as it decays, which is expected to come from the effect of the nonlinearity of the advection term. Experiments on more cases also indicate that the validity of the optimized transmission condition decreases as such discontinuity occurs in their solutions. It noted that, in the method to search for optimal $(\alpha ^*,\beta ^*)$ adopted in this work, simplifications have been adopted, e.g., $\beta ^*=-\alpha ^*$, and this may not result in the real optimal values for $(\alpha, \beta )$. In this sense, it is possible that still the optimized condition remains valid if better approximation for $(\alpha ^*, \beta ^*)$ is adopted. 

\section {Concluding remarks}

This paper studies the computation of ADR equations by the waveform relaxation method. It analyses the convergence speed of the computation of linear equations associated with different transmission conditions, and an optimized transmission algorithm is presented for the Dirichlet and Robin conditions. Then, the optimized algorithm is extended to nonlinear equations. Numerical experiments confirm the analysis and the optimized derived algorithm. This paper leads to the following conclusions. 

1) In the situation of linear equations, the computation converges when the Dirichlet condition is applied. When the Robin condition is adopted, it presents no convergence when the subdomains overlap at two nodes. Nevertheless, it exhibits convergence at speed identical to that associated with the Dirichlet condition with three overlapping nodes. If the combined  condition is adopted, the computation converges faster when the two grids overlap at three nodes than at two nodes.  

2) With the optimized algorithms for the Dirichlet and the Robin conditions, the computation has the same convergence speed if it converges. Numerical experiments show that, indeed, the optimized algorithm speeds up the convergence. 

3) Numerical experiments indicate that the optimized algorithm applies to nonlinear equations, and its performance in speedup largely remains true. The algorithm tends to fail as the effects of the nonlinearity in the convection terms become pronounced. 

Based on this paper, more topics are worthy of further investigation. A topic is to find out why the optimized algorithm fails in nonlinear equations and design schemes to utilize the optimized algorithm for nonlinear equations. Another topic is the extension of this work to higher spatial dimensions, which are commonly the scenarios in real-world problems. A more interesting topic is comparing the waveform relaxation method and the conventional method, which adopts Schwarz iteration between two adjacent time steps, in terms of convergence speed and computational load. We shall keep these as our future research topics. 

\vspace {0.5 cm}

\noindent {\bf Acknowledgments} This work is supported by NSF (DMS $\#$1543876 and $\#$ 1622459).

\vspace {2. cm}

\appendix
\section{Appendix A. Roots of the characteristic equation}

Consider the roots for the characteristic equation, i.e., Eq. (\ref {roots1}). 
%Let $s=\sigma +i\omega $, $\sigma $ and $\omega $ are real numbers, and $\sigma >0$. 
A relation between the two roots is readily derived as 
\begin{equation}
\label{roots4}
r_-r_+=\frac {a}{c}
\end{equation}

\noindent The derivatives of the roots read as  
\begin{equation*}
\label{derivative}
\frac {dr_\pm }{ds}=\-\frac {1}{2c}\mp \frac {s-b}{2c\sqrt {(s-b)^2-4ac}}
\end{equation*}

\noindent Furthermore, one has 
\begin{equation}
\begin{array}{ll}
\label {4acb}
4ac&=4(\dfrac {\mu }{2\Delta x}+\dfrac {\theta }{{\Delta x}^2})
     (-\dfrac {\mu }{2\Delta x}+\dfrac {\theta }{{\Delta x}^2}) \\
   &=\left (\dfrac {2\theta }{{\Delta x}^2}\right )^2-\left ( \dfrac {\mu }{\Delta x} \right )^2 \\
   &=\left (\dfrac {2\theta }{{\Delta x}^2}+\gamma \right )^2-\dfrac {4\theta \gamma }{{\Delta x}^2}-\gamma ^2 -\left ( \dfrac {\mu }{\Delta x} \right )^2 \\
   &<b^2
\end{array}
\end{equation}

\noindent in which $\mu $, $\theta$, $\gamma $=const, and $\theta$, $\gamma $ $>0$. Therefore, it is easy to check that that real part of $\sqrt {(s-b)^2-4ac}$ is positive. As a result, $\sqrt {(s-b)^2-4ac)}\ne 0 $, thus the denominator in the 2nd term of the RHS in the derivatives is not zero. Furthermore, consider $c\ne 0$. All of these mean the derivatives of the roots exit, and thus the roots are analytical \cite {Bak2010}. Since an analytical function takes extremes on its boundary, $r_\pm $ will take minimums and maximums of their values (over the half plane $\sigma \ge 0$) on $\sigma =0$ only, and these minimums and maximums will be analyzed as follows. 

First, let us discuss the roots (\ref {roots1}) in a general situation in which all $\mu $, $\theta $, and $\gamma $ are in presence. In this situation, $b<0$. When $\sigma =0$,  
\begin{equation*}
\begin{array}{ll}
r_\pm &=\dfrac {i\omega -b\pm \sqrt {(i\omega -b)^2-4ac}}{2c} \\
      &=\dfrac {\pm A-b}{2c}+\dfrac {\pm B+\omega }{2c}i
\end{array}
\end{equation*}

\noindent for which $A+Bi=\sqrt {(b-i\omega )^2-4ac} $, and $A=-b\omega /B$. From above, the following is derived:
\begin{equation*}
|r_\pm |^2=\dfrac {A^2+B^2+b^2+\omega ^2\pm 2(B\omega -Ab)}{4c^2}
\end{equation*}

\noindent By some algebra, it can be verified that $A >0$, and $B$ and $\omega$ have a same sign. As a result, one has $B\omega -Ab>0>-(B\omega -Ab)$ and thus $|r_-|^2<|r_+|^2$, or, $|r_-|/|r_+|<1$. Therefore, by (\ref {roots4}), one has 
\begin{equation*}
1>\left |\dfrac {r_-}{r_+}\right |=\dfrac {a}{c}\dfrac {1}{{r_+}^2}
\end{equation*}

\noindent which, in view that $a/c>1$, leads to the second inequality in (\ref {r_-+}). When $\sigma =0$, and under condition (\ref {cond_r_-})  
\begin{equation*}
|r_-|=\dfrac {2a}{|i\omega -b+ \sqrt {(i\omega -b)^2-4ac}|}<1
\end{equation*}

\noindent which is the first inequality in (\ref {r_-+}). 

Note that condition (\ref {cond_r_-}) is satisfied unconditionally in some scenarios. For instance, when $-b>a$, $c<0$, in view of (\ref {abc}), $-b>0$, one has $Re[i\omega -b)^2-4ac]>0$ and  
\begin{equation*}
\begin{array}{ll}
 &|i\omega -b+ \sqrt {(i\omega -b)^2-4ac}| \\
\ge &Re[i\omega -b+ \sqrt {(i\omega -b)^2-4ac}] \\
>&Re[i\omega -b+ \sqrt {(i\omega -b)^2}] \\
\ge &-2b \\
>&2a
\end{array}
\end{equation*}

\noindent Another scenario is $-b>2a$, which leads to 
\begin{equation*}
\begin{array}{ll}
&|i\omega -b+ \sqrt {(i\omega -b)^2-4ac}| \\
\ge &Re[i\omega -b+ \sqrt {(i\omega -b)^2-4ac}] \\
>&Re[i\omega -b] \\
=&-b \\
>&2a
\end{array}
\end{equation*}

\noindent In above discussions, the imaginary parts are excluded, whose increase in value in general reduces the value of $|r_-|$ rapidly. 

In summary, in general, we have $|r_+|>1$, and, under condition (\ref {cond_r_-}), $|r_-|<1 $, along $\sigma =0 $. Since $r_{\pm}$ are analytical and thus their minimums and maximums along $\sigma =0 $ are those for on the right plane, $|r_-|<1 $, $|r_+|>1$, the whole region of $\sigma >0$.

Second, let us discuss two special situations. One of them is $\theta $, $\gamma $=0, and it corresponds to an advection equation. In this situation, $c=-a<0$. Additionally, (\ref {roots4}) leads to $r_-r_+=1$. When $\sigma =0$, 
\begin{equation*}
r_\pm =\frac {i\omega \pm \sqrt {4a^2-\omega ^2}}{2a} 
\end{equation*}

\noindent it is easy to check that still $r_\pm$ are analytical. In case of $4a^2\ge \omega ^2$, one has 
\begin{equation*}
\label{roots5}
|r_\pm |=\sqrt {\frac {\omega ^2+4a^2-\omega ^2}{2a}}=1
\end{equation*}

\noindent In case of $4a^2<\omega ^2$, it yields that 
\begin{equation}
\begin{array}{ll}
|r_+|&=\dfrac {\omega +\sqrt {\omega ^2-4a^2}}{2a}  \\
     &>\dfrac {\omega }{2a} \\
     &>1 
\end{array}
\end{equation}

\noindent As a result, $|r_-|$=$1/|r_+|<1$.

Another special situation is $\mu  =0$, $\gamma $=0, and it corresponds to a diffusion equation. In this situation, $a=c, b=2a$, and (\ref {roots4}) leads to $r_-r_+=1$.  
\begin{equation*}
\label{roots5_2}
r_\pm =\frac {i\omega -b \pm \sqrt {(i\omega -b)^2-4a^2}}{2a} 
\end{equation*}

\noindent Following the derivation in above general situation, it can be shown that $|r_+|<1$. Then, we have $|r_-|$=$1/|r_+|<1$.  

\let\oldbibliography\thebibliography
\renewcommand{\thebibliography}[1]{\oldbibliography{#1}
\setlength{\itemsep}{-3pt}}
\bibliography{references}
\bibliographystyle{ieeetr}

\end{document}